\documentclass[11pt]{article}
\usepackage{latexsym,amsmath,amssymb,amscd, graphics}
\usepackage{pslatex}
\usepackage{amsmath}

\makeatletter \@addtoreset{figure}{section}
\def\thefigure{\thesection.\@arabic\c@figure}
\def\fps@figure{h,t}
\@addtoreset{table}{bsection}

\def\thetable{\thesection.\@arabic\c@table}
\def\fps@table{h, t}
\@addtoreset{equation}{section}

\makeatother

\newtheorem{theorem}{Theorem}

\newtheorem{corollary}[theorem]{Corollary}

\newtheorem{definition}[theorem]{Definition}

\newtheorem{lemma}[theorem]{Lemma}

\newtheorem{proposition}[theorem]{Proposition}
\newtheorem{remark}[theorem]{Remark}

\numberwithin{theorem}{section} 
\begin{document}
\title{Metric convexity in the symplectic category}
\author{Petre Birtea, Juan-Pablo Ortega, Tudor S. Ratiu}
\date{September 18,   2006}
\maketitle

\begin{abstract}
We introduce an extension of the standard Local-to-Global Principle used in the proof of the
convexity theorems for the momentum map to handle closed maps that take values in a length
metric space. This extension is used to study the convexity properties of the cylinder
valued momentum map introduced by Condevaux, Dazord, and Molino in~\cite{condevaux dazord
and molino} and allows us to obtain the most general convexity statement available in the
literature for momentum maps associated to a symplectic Lie group action. 
\end{abstract}

\medskip


\medskip

\section{Introduction}

The study of the convexity properties of the image of momentum maps has been a very active
field of research for the last twenty years ever since Atiyah \cite{atiyah 82}
and, independently, Guillemin and Sternberg \cite{convexity}, realized that this is a natural
way to encode very important problems in linear algebra and representation theory. One of
them is the classical result of Schur \cite{schur} and
Horn \cite{horn} that states that the set of diagonals of an
isospectral set of $n\times n$ Hermitian matrices equals the convex hull of
the $n!$ points obtained by permuting all the eigenvalues. 
A second linear algebra problem that can be described in the momentum map setup is the
characterization of all possible eigenvalues of the sum 
$A + B $ of two Hermitian matrices $A$ and $B
$ as each one of them ranges over an isospectral set. The isospectral sets
of Hermitian matrices are precisely the coadjoint orbits
$\mathcal{O}_\mu$ of $U(n)$, where $\mu = (\mu_1, \mu_2, \dots,
\mu_n)$ and $\mu_i $ are the eigenvalues.  If one requires, in
addition, that the eigenvalues of the sum $A + B $ be sorted in
decreasing order, this problem amounts to describing the intersection of the image of a
coadjoint equivariant  momentum map $\mathbf{J}:M \rightarrow \frak{u}(n) ^\ast  $ with one
of the Weyl chambers of
$\frak{u}(n) ^\ast$. More explicitly, in this case $M$ is the symplectic manifold obtained by
taking the Cartesian product of the two $U (n) $ coadjoint orbits $\mathcal{O}_{\mu}\times
\mathcal{O}_{\lambda } $ and  $\mathbf{J}:M \rightarrow \frak{u}(n) ^\ast  $ is the
momentum map corresponding to the diagonal $U(n) $-action on that set.

This
problem is a particular case of a more general
situation  which consists of describing the image of a coadjoint equivariant momentum map
$\mathbf{J}:M\rightarrow \mathfrak{g}^{*}$ associated to the action of a compact Lie group
$G$  on a compact symplectic manifold $M$.
Guillemin and Sternberg \cite{convexity2} proved that $\mathbf{J}(M)\cap\mathfrak{t}^{*}_+$
is a union of compact convex polytopes ($\mathfrak{t} $ is a maximal toral algebra of
$\mathfrak{g}$ and   $\mathfrak{t}^{*}_+ $ the positive Weyl chamber of
$\mathfrak{g}^\ast $) and Kirwan
\cite{kirwan convexity} showed that this set is connected thereby concluding that
$\mathbf{J}(M)\cap
\mathfrak{t}_{+}^{*}$ is a compact convex polytope.

These convexity results have been generalized to compact group actions on noncompact
manifolds with proper momentum maps by Condevaux, Dazord, and Molino \cite{condevaux
dazord and molino} who designed a very ingenious topological proof of the 
convexity theorems based on the local properties of the
momentum map. This method is usually referred to as Local-to-Global
Principle. The main tool for capturing the local properties of the momentum map was the
so-called  Marle-Guillemin-Sternberg normal form
\cite{marle, GuSt1984}. Using this approach they proved that a momentum map
of a torus Hamiltonian action is locally fiber connected, locally open
onto its image, and has local convexity data. This method was
further analyzed and applied to many interesting situations by
Hilgert, Neeb, and Plank \cite{hilnebplank}. Later on Sjamaar
\cite{sjamaar} and  Knop \cite{knop} proved that a momentum
map in the non-Abelian case has the same local topological
properties mentioned above.

In spite of its generality, the Local-to-Global Principle cannot be
applied to obtain a convexity result which would include other interesting
examples like, for example, those introduced by Prato \cite{prato} since the momentum map in
this case is not proper. This inconvenience has been fixed by the authors in~\cite{polytope
paper} where the properness condition has been replaced  by just closedness.

All the above convexity results concern maps with values in
vector spaces. Consequently, they address the situation in which
the group action in question has an associated standard momentum map. 
When such a map does not exist one still can define a momentum map that captures most of the
properties of the standard object but which, in general, takes values in an Abelian group
isomorphic to a cylinder. 
This momentum map that will be referred to as the \textbf{cylinder valued momentum map}, was
introduced in Condevaux, Dazord, and Molino \cite{condevaux dazord and molino}
and carefully studied in Ortega and Ratiu \cite{cylinder reduction, stratified cylinder
valued momentum map} in the context of reduction. Additionally, its local properties are as
well known as those for the standard momentum map~\cite{symplectic slice, stratified
cylinder valued momentum map}. Hence, it is very natural to ask if one could extend
to the cylinder valued momentum map the knowledge that we have about the convexity properties
of the classical momentum map by using an appropriate generalization of the Local-to-Global
Principle. 
This question is actually posed as an open problem in the original article \cite{condevaux dazord and 
molino}.
An affirmative answer to this twenty year old question is the main achievement of this paper.

The main results of the present work are divided into two parts. First, 
in Section~\ref{Image convexity for maps with values in length spaces} we extend the Local-to-Global 
Principle to the
category of maps that take values in length spaces. Our approach  is based on the
extension  of the classical Hopf-Rinow theorem to length metric spaces due to Cohn-Vossen. This circle of ideas
also appears for compact spaces in~\cite{buliga}. The  results on length metric spaces necessary for the comprehension of the paper are
recalled in an appendix at the end.  The metric approach seems to be the best adapted
generalization of the classical setup to our problem since, under certain hypotheses related
to the topological nature of the Hamiltonian holonomy of the problem (a concept  defined carefully later 
on), the target space of the cylinder valued momentum map has an associated
canonical length space structure.  The extension presented here recovers the previous Local-to-Global Principle proved
in~\cite{polytope paper} by the authors for the case of closed maps with values in vector spaces. It is 
expected that this general result can be applied to many other situations going beyond symplectic 
geometry. Second, in Section~\ref{Metric convexity for cylinder valued momentum maps} we  apply
this generalized Local-to-Global Principle to the  cylinder valued momentum map and we
obtain a convexity result similar to the classical one for the standard momentum map.

\section{Image convexity for maps with values in length spaces}
\label{Image convexity for maps with values in length spaces}

One of the main goals in this paper is the study of the convexity properties of the image of a natural generalization of the momentum map. The notion of convexity is usually associated with vector spaces. However, the map considered in this paper has values in a manifold that is,
in general, diffeomorphic to a cylinder. Thus, one  is forced to work in a more general setting. As reviewed in the appendix 
\S\ref{Appendix: metric and length spaces}, most of the concepts pertaining to convexity can be extended to the context of the so called 
\textbf{length spaces}. It turns out that the target space of the map that we are going to study can be naturally endowed with a length space structure and hence convexity will be used in this context. We give in 
\S\ref{Appendix: metric and length spaces} a self-contained brief summary of  all the definitions and results on length spaces necessary in this paper.

The convexity program has been successfully carried out 
for the standard momentum map by several means. One possible approach consists in determining certain local properties of the map  that guarantee that it has a globally convex image. This strategy relies on a fundamental result called the \textbf{Local-to-Global Principle} which has been introduced in~\cite{condevaux dazord and molino, hilnebplank} for maps whose target space is an Euclidean vector space. Since the extension of the standard momentum map with which we
will be working does not map into a vector space but into a length space, a
generalization of the Local-to-Global
Principle is needed to handle this situation. This is the main goal of the present section.
\medskip

Let $f:X\rightarrow Y$ be a continuous map between two connected
Hausdorff topological spaces.  Define
the following equivalence relation on the topological space $X$: 
declare two points $x_1,x_2\in
X$ to be equivalent if and only if $f(x_1)=f(x_2)=y$ and they
belong to the same connected component of $f^{-1}(y)$. The
topological quotient space, whose elements are the connected components 
of the fibers of $f $, will be denoted by $X_{f}$, the
projection map by $\pi _{f}:X\rightarrow X_{f}$, and the induced
map on $X_{f}$ by $\widetilde{f}:X_{f}\rightarrow Y$. Thus,
$\widetilde{f} \circ \pi_f = f$ uniquely characterizes
$\widetilde{f}$. The map $\widetilde{f}$ is continuous and if the
fibers of $f$ are connected then it is also injective.

\begin{definition}\label{LFC}
Let $X $  and $Y$  be two  topological spaces and $f:X \rightarrow
Y$ a continuous map. The subset $A \subset X$ satisfies the 
\textbf{locally fiber connected condition (LFC)} if $A$ does not intersect
two different connected components of the fiber $f ^{-1}(f (x))$,
for any $x \in A$.

Let $X $ be an arcwise connected Hausdorff topological
space. The continuous map $f: X \rightarrow Y$ is said to be
\textbf{locally fiber connected} if for  each $x\in X$, any open
neighborhood of $x$ contains a connected neighborhood $U_{x}$ of
$x$ such that $U_{x}$ satisfies the \textbf{(LFC)} condition.
\end{definition}

The following consequences of the definition will be useful later on.
A subset of a set that satisfies (LFC) also satisfies (LFC). If $A
\subset X $ satisfies the (LFC) property, then its saturation
$\pi_f^{-1}(\pi_f(A))$ also satisfies (LFC). If $f $ is locally fiber
connected, then any open neighborhood  of $x\in X$  contains an open
neighborhood $U_x $  of
$x$ such that the restriction of $\widetilde{f} $  to $\pi _f(U _x)$ 
is injective. 

\begin{definition}
A continuous map $f: X \rightarrow Y$ is said to be
\textbf{locally open onto its image} if for any $x\in X$ there
exists an open neighborhood $U_x$ of $x$ such that the restriction
$f|_{U_x}:U_x \to f(U_x)$ is an open map, where $f(U_x)$ has 
the  topology induced by $Y $. We say that such a
neighborhood satisfies the \textbf{(LOI)} condition.
\end{definition}

\noindent Benoist proved in Lemma 3.7 of~\cite{benoist} the following result that will be used
later on.

\begin{lemma}\label{benoist}
Suppose $f:X\rightarrow Y$ is a continuous map between two
topological spaces. If $f$ is locally fiber connected and locally
open onto its image then $\pi_f$ is an open map.
\end{lemma}

\noindent The following characterization of closed maps will be needed in what follows (see \cite{engelking}, Theorems 1.4.12 and 1.4.13).

\begin{theorem}
\label{caracterizare inchise} Let $f:X \rightarrow Y $ be a
continuous mapping.
\begin{description}
\item[{\rm \textbf{(i)}}] $f$ is closed if and only if for every
$B\subset Y$ and every open set $A\subset X$ which contains
$f^{-1}(B)$, there exists an open set $C\subset Y$ containing $B$
and such that $f^{-1}(C)\subset A$.

\item[{\rm \textbf{(ii)}}] $f$ is closed if and only if for every
point $y\in Y$ and every open set $U\subset X$ which contains
$f^{-1}(y)$, there exists a neighborhood $V_y$ of the point
$y$ in $Y$  such that $f^{-1}(V_y)\subset U$.
\end{description}
\end{theorem}

For the next results we recall the following standard definitions
(see, e.g., Engelking \cite{engelking}).  A topological space $X$ is
called a
$\operatorname{T}_1$-space if for every pair of distinct points $x
_1, x _2 \in X $, there exists an open set $U \subset X $ such that
$x _1 \in U $ and $x _2\notin U $. The topological space $X $ is
\textbf{normal} if it is a $\operatorname{T} _1$-space and for any
closed disjoint subsets $A, B \subset X$ there exist open subsets $U,
V \subset X$ such that $A \subset U$, $B \subset V $, and $U \cap V =
\varnothing$. The topological space $X $ is \textbf{first countable}
if every point admits a countable base of open neighborhoods.

\begin{lemma}
\label{for Hausdorff} 
Let  $X$ be a normal, first countable, arcwise connected, and  Hausdorff
topological space and $Y$ a Hausdorff topological space. Let
$f:X\rightarrow Y$ be a continuous map that is locally  open onto
its image and is locally fiber connected. If $f$ is a closed
map, then
\begin{description}
\item[(i)] the projection $\pi_f: X \rightarrow X _f$ is also a
closed map,

\item[(ii)] the quotient  $X_f$ is a Hausdorff topological space.
\end{description}
\end{lemma}

\noindent\textbf{Proof.\ \ } \textbf{(i)} Let $[x]$ be an arbitrary point
in $X_f$ and $U\subset X$ an arbitrary open set that includes
$E_x: = {\pi_f}^{-1}([x])$, the connected component of $f^{-1}(f(x))$
that contains $x$. Denote by $F:=f^{-1}(f(x))\setminus E_x$ the union of all (closed) connected components of
$f^{-1}(f(x))$ different from $E_x$. We claim that $F$ is a closed subset
of $X$. Indeed, if $z\in \overline{F}$, by first countability of
$X $, there exists a sequence $\{z_n\}_{n \in  \Bbb N}$ in $F$ which
is convergent to $z$. Since $f(z _n) = f(x)$, by continuity of $f$ we
conclude that $f (x) = f( z _n) \to f(z)$ and hence
$z\in f^{-1}(f(x))$. If $z\in E_x$ then any neighborhood of $z$
intersects at least one other connected component of the fiber
$f^{-1}(f(x))$ since $z \in \overline{F}$. This, however,
contradicts the (LFC) condition. Therefore, $ z \in F$ 
and hence $F$ is closed. The same argument as above shows that the (LFC)
condition implies that
$E_{x}$ is also closed in $X$.

Using the normality of $X$ there exist two open sets $U_{E_x}$ and
$W$ such that $E_x\subset U_{E_x}$, $F\subset W$, and
$U_{E_x}\cap W=\varnothing$. After shrinking, if necessary, we can
assume that $U_{E_x}\subset U$. Applying 
Theorem~\ref{caracterizare inchise}\textbf{(ii)}, the closedness of
$f$ ensures the existence of an
open neighborhood  $V_{f (x)}$ of  $f (x)  $ in $V$ such that
$E_x\subset f^{-1}(f(x))\subset f^{-1}(V_{f (x)})\subset
U_{E_x}\cup W$.

The set $A:=U_{E_x}\cap f^{-1}(V_{f (x)})$ is a nonempty open
subset of $X$ and is also saturated with respect to the
equivalence relation that defines $\pi_f$ or, equivalently,
${\pi_f}^{-1}(\pi_f(A))=A$. Indeed, if a connected component of a
fiber of $f$ from $f^{-1}(V_{f (x)})$ intersects $U_{E_x}$,
respectively $W$, then it is entirely contained either in
$U_{E_x}$ or in $W$ since $U_{E_x}\cap W=\varnothing$. 

Since $A $ is open in $X $, by the definition of the quotient
topology of $X _f$, it follows that $\pi_f(A)$ is an open
neighborhood of $[x]$. Note that ${\pi_f}^{-1}(\pi_f(A))\subset U$,
which shows via Theorem \ref {caracterizare inchise}\textbf{(ii)}
that $\pi_f$ is a closed map.

\medskip

\textbf{(ii)} We shall prove that $X _f$ is Hausdorff by showing  that the
projection $\pi _f$ is an open map (which holds by Lemma
\ref{benoist}) and that the graph of the equivalence relation that
defines $X _f$ is closed.

To show that the graph is closed, we need some preliminary
considerations. For every $x\in X$ there exists a
neighborhood $U_x$ that satisfies (LOI) and (LFC). By normality of $X$
there exists also a neighborhood $U'_x$ of $x$ with
$\overline{U'_x}\subset U_x$. We shall prove that
$\overline{{\pi_f}^{-1}(\pi_f(U'_x))}\subset
{\pi_f}^{-1}(\pi_f(U_x))$ which shows that for every connected
component $E_x$ of a fiber there exists a saturated neighborhood of it
which contains a smaller saturated neighborhood whose closure still
satisfies (LFC). In order to prove the above inclusion observe that
since $\pi_f$ is continuous and closed we have that
$\overline {\pi_f(U'_x)}=\pi_f(\overline{U'_x})\subset \pi_f(U_x)$.
By the continuity of $\pi_f$ we obtain the inclusion
$\overline{{\pi_f}^{-1}(\pi_f(U'_x))}\subset
{\pi_f}^{-1}(\overline{\pi_f(U'_x)})\subset
{\pi_f}^{-1}(\pi_f(U_x))$.

We now prove the closedness of the graph of the equivalence
relation that defines $X _f  $. Take $\{x_n\}_{n\in\mathbb{N}}$
and $\{y_n\}_{n\in \mathbb{N}}$ two convergent sequences in $X$
such that $x_n$ and $y_n$ are in the same equivalence class for
all $n\in \mathbb{N}$. Suppose that $x_n\rightarrow x$ and
$y_n\rightarrow y$. The continuity of $f$ guarantees that $f  (x)
= f  (y) $. Additionally, there exists $n_0\in \mathbb{N}$ such
that for $n>n_0$ all $x_n\in {\pi_f}^{-1}(\pi_f(U'_x))$, where $U
_x'$ has been chosen as above. Consequently $y_n\in
{\pi_f}^{-1}(\pi_f(U'_x))$ since $x_n$ and $y_n$ are in the same
equivalence class and ${\pi_f}^{-1}(\pi_f(U'_x))$ is saturated.
Therefore, $x, y \in \overline{{\pi_f}^{-1}(\pi_f(U'_x))}$. But
$\overline{{\pi_f}^{-1}(\pi_f(U'_x))}$ satisfies (LFC) and thus
$x $ and $y $ belong to the same connected component of the fiber
$f^{-1}( f(x))$. This shows that the graph of the equivalence
relation is closed, as required.
\quad $\blacksquare$

\medskip

On $X_f$ we define the function $\widetilde d:X_f\times
X_f\rightarrow [0, \infty]$ in the following way: for $[x]$ and $[y]$
in $X_f$ let $\widetilde d ([x], [y])$ be the infimum of all the
lengths $l_d (\widetilde f \circ \gamma)$ where $\gamma$ is a
continuous curve in $X_f$ that connects $[x]$ and $[y]$. The
length $l_d$ is computed with respect to the distance $d$ on
$Y$. From  the definition it follows that $d(\widetilde f([x]),
\widetilde f ([y]))\leq \widetilde d ([x], [y])$.

\begin{proposition}
Let $X$ be a normal, first countable, arcwise connected, and  Hausdorff topological space
and $(Y, d)$ a metric space. Assume that
$f:X\rightarrow Y$ is a continuous closed map that is also locally open onto its image and
locally fiber connected. Then the function on  $\widetilde d:X_f\times
X_f\rightarrow [0, \infty]$ introduced
above defines a metric on $X_f$.
\end{proposition}

\noindent\textbf{Proof.\ \ }The positivity, symmetry, and the triangle inequality of
$\widetilde d$ are obvious from the definition of $\widetilde d$. It remains to be
proved that if $\widetilde d ([x], [y]) =0$ then $[x]=[y]$.
Suppose that there exist $[x] \not = [y]$ with $\widetilde d ([x],
[y])=0$. Then $d(\widetilde f ([x]), \widetilde f([y]))=0$ and
hence $f(x)  = \widetilde f([x])= \widetilde f([y]) = f (y)$. This
implies that
$[x]$ and $[y]$ are images under the the projection $\pi_f$ of two
different connected components of the same fiber.

By the (LOI) property of $f$ and the openness of $\pi_f$ (see
Lemma \ref{benoist}) there exist two open neighborhoods $U_{[x]}$ and
$U_{[y]}$ in the quotient topology of $X_f$ around $[x]$,
respectively $[y]$, such that $\widetilde{f}|_{U_{[x]}}$ and
$\widetilde{f}|_{U_{[y]}}$ are open onto their images. Consequently,
there exist two open neighborhoods $U^\prime_{[x]}\subset U_{[x]}$
and $U^\prime_{[y]} \subset U_{[y]}$ such that $\widetilde{f}
(U^\prime_{[x]})\supset B_d (\widetilde f[x], r) \cap f (U_{[x]})$ and
$\widetilde{f}(U^\prime_{[y]})\supset B_d (\widetilde f
([y]),r^\prime)\cap f(U_{[y]})$ for two small enough constants $r$,
$r^\prime>0$ . Moreover, we can choose the above
neighborhoods $U^\prime_{[x]}$ and $U^\prime_{[y]}$ such that
$U^\prime_{[x]}\cap U^\prime_{[y]}=\varnothing$ since the quotient topology of $X_f$ is Hausdorff by Lemma \ref{for Hausdorff}.

Any curve $\gamma$ in $X_f$ that connects $[x]$ and $[y]$ is
mapped by $\widetilde {f}$ into a loop in $Y$ based at $\widetilde
f([x])=\widetilde f([y])$ and has to exit $U^\prime_{[x]}$ and
enter $U^\prime_{[y]}$ (since $U^\prime_{[x]}\cap
U^\prime_{[y]}=\varnothing$) in order to connect $[x]$ and $[y]$. Consequently, the curve $\widetilde{f} \circ\gamma$ in $Y$ has to exit $B_d (f([x]),r)$ and reenter
$B_d(f([y]),r^\prime)$ and hence $l_d(\widetilde{f} \circ \gamma)>
r+r^\prime$.
This is in contradiction with the hypothesis that
$\widetilde d([x], [y])=0$ for $[x] \not = [y]$.
\quad $\blacksquare$

\medskip 

In order to put the following definition in context, the reader is encouraged to look at the appendix \S\ref{Appendix: metric and length spaces} where the concepts of length metric and geodesic metric space are discussed.

\begin{definition}\label{metric convex 1}
A subset $C$ in a metric space $(X, d)$ is said to be \textbf{convex}
if the restriction of $d$ to $C$ is a finite length metric 
{\rm (}see Definition {\rm \ref{length space})}.
\end{definition}

\begin{definition}\label{local convex} Let $X$ be a connected Hausdorff
space and $(Y,d)$ a length space {\rm (}see Definition {\rm \ref{length space})}. A continuous mapping $f:X \to Y$
is said to have \textbf{local convexity data} if for each $x \in X$
and every sufficiently small neighborhood $U_x$ of $x$ the set
$f(U_x)$ is a convex subset of $Y$. We say that $U_x$ satisfies the
\textbf{(LCD)} condition.
\end{definition}

\begin{proposition}\label{equality}
Let $X$ be a normal, first countable, arcwise connected, and  Hausdorff
topological space and $(Y,d)$ a geodesic metric space {\rm (}see Definition {\rm \ref{geodesic})}. Assume that
$f:X\rightarrow Y$ is a continuous closed map that is also locally
open onto its image, locally fiber connected, and has local
convexity data. Then $(X_f, \widetilde d)$ is a length space and
the topology induced by $\widetilde d$ coincides with the quotient
topology of $X_f$.
\end{proposition}

\noindent\textbf{Proof.\ \ } First we will prove that for close enough $[x]$ and
$[y]$ we have $d(\widetilde f([x]), \widetilde f ([y])) =
\widetilde d ([x], [y])$. Indeed, let
$[x], [y] \in X_f$ be such that they are contained in an open
set $\widetilde{U} \subset X_f$ (open in the quotient topology of
$X_f$) such that $\widetilde{U} = \pi _f(U)$, where $U$ open in $X $ and
satisfies the (LCD), (LOI), and (LFC) conditions. Let $\Omega : =
f(U) = \widetilde{f}( \widetilde{U})$. Since
$\Omega$ is convex (because (LCD) holds for $U $), by Lemma
\ref{caracterizare} there exists a rectifiable shortest path $\gamma_0 $
entirely contained in $\Omega$ that connects $\widetilde f([x])$
and $\widetilde f([y])$, that is,
$l_d(\gamma_0)=d(\widetilde{f}([x]),\widetilde{f}([y]))$.   Note that
$\widetilde f|_{\widetilde{U}} :\widetilde{U} \rightarrow \Omega$ 
is a homeomorphism ($\widetilde{U}$ endowed with the quotient topology
of $X_f$) because $\widetilde f|_{\widetilde{U}}$ is open, since  $U $
satisfies (LOI), and is injective, because $U $ satisfies (LFC). The
curve $ c_0:=\widetilde f^{-1} (\gamma_0)$  is continuous and
connects $[x]$ with $[y]$. From the definition of $\widetilde{d}$ we
have that $\widetilde{d}([x],[y])\leq l_d(\widetilde{f}\circ
c_0)=l_d(\gamma_0)=d(\widetilde{f}([x]),\widetilde{f}([y]))$. As
the inequality $d(\widetilde f([x]), \widetilde f ([y]))\leq
\widetilde d ([x], [y])$ always holds,  we obtained the desired
equality $d(\widetilde f([x]), \widetilde f ([y])) =
\widetilde d ([x], [y])$. Consequently, the homeomorphism
$\widetilde f|_{\widetilde{U}}: \widetilde{U} \to \Omega$ is also an
isometry from $(\widetilde{U}, \widetilde d|_{\widetilde{U}})$ to
$(\Omega, d|_{\Omega})$ and thus the quotient topology of $X_f$
coincides with the metric topology induced by $\widetilde d$.

Next we will prove that $(X_f, \widetilde d)$ is a  length space.
Let $c:[a,b]\to X_f$ be a continuous curve connecting two
arbitrary points $[x]$ and $[y]$ in $X_f$. For two partitions
$\Delta_n$ and $\Delta_{n+1}$ of the interval $[a,b]$ with
$\Delta_{n+1}$ finer then $\Delta_{n}$ we have that $\sum_{i=1}^n
\widetilde d (c(t_i), c(t_{i+1}))\leq \sum_{i=1}^{n+1} \widetilde d
(c(t_i), c(t_{i+1}))$ due to the triangle inequality.  Therefore, 
in order to compute $l_{\widetilde d}(c)$ it suffices to work with
partitions fine enough such that two consecutive points $c(t_i),
c(t_{i+1})$, corresponding to a partition $\Delta_n$, are close enough 
as above. Therefore, by what was just proved, we have
$d (\widetilde{f}(c(t_i)), \widetilde f(c(t_{i+1}))) = \widetilde{d}
(c(t_i), c(t_{i+1}))$ and we conclude
\begin{eqnarray*}
l_{\widetilde d} (c) &=& \sup_{\Delta_n} \sum_1^n \widetilde
d (c(t_i), c(t_{i+1})) \\
&=& \sup_{\Delta_n}\sum_1^n d (\widetilde f
(c(t_i)), \widetilde f(c(t_{i+1}))) \\
&=& l_d (\widetilde f \circ c).
\end{eqnarray*}
Consequently, $\widetilde f \circ c$ is a rectifiable curve in
$(Y, d)$ if and only if $c$ is rectifiable in
$(X_{f},\widetilde{d})$. The equality
$$
\widetilde d ([x], [y]) = \inf l_d (\widetilde f \circ c) = \inf
l_{\widetilde d} (c) = \overline{\widetilde d}([x], [y]),
$$ shows
that $(X_f, \widetilde d)$ is a length space.
\quad $\blacksquare$

\medskip

The proof of the following lemma can be found as Proposition 4.4.16 in~\cite{engelking}.

\begin{lemma}
\label{frontier} (Va\v{\i}n\v{s}te\v{\i}n) If $f:X\rightarrow Y$
is a closed mapping from a metrizable space $X$ onto a metrizable
space $Y$, then for every $y\in Y$ the boundary ${\rm bd
}(f^{-1}(y)):=\overline{f^{-1}(y)}\cap\overline{ \left(X \setminus
f^{-1}(y)\right)}$ is compact.
\end{lemma}

\begin{definition}\label{proper map} Let $X$ be a Hausdorff
topological space and $f:X \to Y$ a continuous map. We call $f$
a \textbf{proper map} if it is closed and all fibers $f^{-1}(y)$ are
compact subsets of $X$.
\end{definition}

The proof of the following theorem can be found in Engelking
\cite{engelking}, Proposition 3.7.2.

\begin{theorem}
\label{proper}
If $f:X \to Y$ is a proper map,
then for every compact subset $Z\subset Y$ the inverse image
$f^{-1}(Z)$ is compact.
\end{theorem}

A converse of this theorem is available when $Y$ is a $k$-space
(i.e $Y$ is a Hausdorff topological space that is the image of a
locally compact space under a quotient mapping). For example every
first countable Hausdorff space is a $k$-space (see Theorem 3.3.20 of~ \cite{engelking}).

\begin{proposition}
Let $X $ be a normal, first countable, arcwise connected, and  Hausdorff topological space
and $(Y, d)$ a complete locally  compact length space (and thus, by Hopf-Rinow-Cohn-Vossen a
geodesic metric space; see Theorem \ref{Hopf-Rinow}). Assume that $f: X \rightarrow Y $ is a continuous closed map that is
also locally open onto its image, locally  fiber connected, and has local convexity data.
Then
$(X, \widetilde{d})$ is a complete locally compact  length space and hence a geodesic
metric space.
\end{proposition}

\noindent\textbf{Proof.\ \ }First we will prove that $\widetilde f$ is a proper map. It is a
closed map since $f$ is a closed map. Local injectivity of
$\widetilde f$ implies that ${\rm bd }(\widetilde
{f}^{-1}(y))=\widetilde {f}^{-1}(y)$. By Va\v{\i}n\v{s}te\v{\i}n's
Lemma \ref{frontier} we conclude that the fibers of $\widetilde{f}$ are all compact
and consequently $\widetilde f$ is a proper map.

Because local compactness is an inverse invariant for proper maps
we obtain also that $(X_f, \widetilde d)$ is locally compact since
$\widetilde f : X_f \to Y$ is a proper map.

By the Hopf-Rinow-Cohn-Vossen Theorem \ref{Hopf-Rinow} it suffices to show that every
closed metric ball in $X_f$ is compact in order to conclude  that
$(X_f, \widetilde d)$ is a complete metric  space. The set
$\overline B ([x], r)$ is closed in $X_f$ and, by definition of the
metric
$\widetilde d$, the inclusion $\widetilde f (\overline B ([x], r))
\subset \overline B (\widetilde f ([x]), r)$ holds. By the
Hopf-Rinow-Cohn-Vossen Theorem it follows that 
$\overline B(\widetilde f ([x]),r)$ is a compact set in $Y$ and, consequently, $\widetilde f^{-1}(\overline B (\widetilde f ([x]), r)$ is compact in $X_f$ due to properness of $\widetilde f$. Since $\overline B([x], r)$ is a
closed subset of $\widetilde f^{-1} (\overline B (\widetilde f
([x])), r)$ it is necessarily compact in $X_f$.
\quad $\blacksquare$

\medskip

As a consequence of the above proposition, $(X_f, \widetilde d)$
satisfies all the conditions of the Hopf-Rinow-Cohn-Vossen Theorem
which implies that for any two points $[x], [y] \in X_f$ there
exists a shortest geodesic connecting them.

\begin{definition}
Let $(X, d)$ be a geodesic metric space. We say that
$C$ is \textbf{weakly convex} if for any two points $x,y\in C$ there
exist a geodesic between $x$ and $y$ entirely contained in $C$.
\end{definition}

Note that weak convexity does not require that the geodesic
to be a shortest one. Now we can present the main result of this
section.

\begin{theorem}[Local-to-Global Principle]
\label{lokal to global metric}
Let $X$ be a normal, first countable, arcwise connected, and  Hausdorff
topological space and $(Y,d)$ a complete, locally compact
length space. Assume that $f:X\rightarrow Y$ is a continuous
closed map that is also locally open onto its image, locally fiber
connected, and has local convexity data. Then the following hold:
\begin{description}
  \item [(i)] $f(X) \subset (Y,d)$ is a weakly convex subset of $Y$.
  \item [(ii)] If, in addition,  $(Y,d)$ is uniquely geodesic (that
is, any two points can be joined by a unique geodesic) and $d(y _1, y
_2) < \infty$ for all $y _1, y _2 \in Y$, then
$f(X)$ is a convex subset of $(Y,d)$, $f$ has connected fibers, and
  it is open onto its image.
\end{description}
\end{theorem}

\noindent\textbf{Proof.\ \ } \textbf{(i)} We have to prove that for any $y_1, y_2 \in f(X)$
there exists a geodesic (not necessary shortest) in $(Y, d)$ completely
included in $f(X)$. Indeed, take $[x_1], [x_2]\in X_f$ such that
$\widetilde f ([x_1]) = y_1$ and $\widetilde f([x_2]) = y_2$. As
was explained above, there exists a shortest geodesic $c:[a,b] \to
X_f$ with the properties $c(a) = [x_1]$, $c(b) = [x_2]$, and
$\widetilde d ([x_1], [x_2]) = l_{\widetilde d}(c)$. We will show
that $\widetilde f \circ c \subset f(X)$ is a geodesic that
connects $y_1$ with $y_2$. Since $f $ is locally open onto its image, 
locally fiber connected, and has local convexity data, each $[x] \in
X _f$ admits an open neighborhood $U_{[x]}$ such that
$\widetilde{f}|_{U_{[x]}}: U_{[x]} \rightarrow
\widetilde{f}(U_{[x]})$ is injective, open onto its image, and
$\widetilde{f}(U_{[x]})$ is convex in $Y $. Choose now $[x]$ in
the image of $c $ and let $t_0$ be such that $c(t_0) = [x]$. Then the
intersection of the image of $c $ with $U_{[x]}$ is the image of a
curve $c^\prime : I \to U_{[x]}$
with $I \subset [a,b]$ a subinterval. If we take $U_{[x]}$ small
enough then for any subinterval $[t_1, t_2] \subset I$ with $t_0\in
[t_1, t_2]$ we have $\widetilde d (c(t_1), c(t_2))= d(\widetilde f (c(t_1)),
\widetilde f(c(t_2)))$ as was explained in the proof of
Proposition \ref{equality}. Because $c$ is a shortest geodesic we have
that $\widetilde d (c(t_1), c(t_2))= l_{\widetilde d}(c|_{[t_1, t_2]})$
and since $l_{\widetilde d}(c|_{[t_1, t_2]})= l_d (\widetilde f
\circ c|_{[t_1, t_2]})$ (as was shown in the proof of Proposition
\ref{equality}) we obtain the desired equality
$d(\widetilde f (c(t_1)), \widetilde f(c(t_2)))=l_d (\widetilde f
\circ c|_{[t_1, t_2]}) $. This proves that $\widetilde f \circ
c$ is a geodesic in $(Y, d)$ since it is a local distance minimizer.

\medskip

\textbf{(ii)} From (i) we already know that any two points in $f(X)$ can be
connected by a geodesic included in $f(X)$. Since $Y$  is uniquely
geodesic this chosen geodesic must be the shortest one. Also, the 
restriction of $d$ to $f(X)$ is finite and therefore
$f(X)$ is a convex subset of $(Y,d)$.

Now we will prove that $f$ has connected fibers. Suppose the
contrary, that is, there exist $[x] \neq [y]$ with $\widetilde f
([x]) =
\widetilde f([y])$. Then there exists a shortest geodesic
$c:[a,b]\rightarrow X_f$ that links $[x]$ and $[y]$ which is
mapped by $\widetilde f$ to a loop based at $\widetilde f([x]) =
\widetilde f([y])$. As was proved in (i), $\widetilde f \circ c$ is
a geodesic in $Y$. Since $(Y,d)$ is uniquely geodesic we obtain
that $\widetilde f\circ c$ is the constant loop. Consequently,
$\widetilde f(c(t))=\widetilde f([x])$ for all $t \in [a,b]$. This
implies that $c(t)$ and $[x]$ belong to the same fiber of
$\widetilde f$  for all $t \in [a,b]$ which contradicts the local
injectivity of $\widetilde f$ implied by the (LFC) property of $f$.

Since $\widetilde f$ is a closed injective map, it is also open
onto its image and, consequently, $f$ is open onto its image.
\quad $\blacksquare$

\begin{remark}
\normalfont
Unlike the situation encountered in the classical Local-to-Global
Principle~\cite{condevaux
dazord and molino, hilnebplank} in which the target space of the map is a Euclidean vector space and hence
uniquely geodesic, $f$ could have, in general, a weakly convex image but not connected
fibers. See Section~\ref{Convexity for Abelian Lie group valued momentum maps} for an example. 
\end{remark}

\begin{remark}
\label{restriction}
\normalfont 
If $Y$ is an Euclidean vector space and $C$ is a convex subset of $Y$ then 
Theorem \ref{lokal to global metric} applied to the map $f:X \rightarrow C$ yields 
the generalization of the classical Local-to-Global Principle introduced in Theorem 2.28 of~\cite{polytope paper}.
\end{remark}

\section{Metric convexity for cylinder valued momentum maps}
\label{Metric convexity for cylinder valued momentum maps}

The goal of this section is to apply the general results obtained in 
\S\ref{Image convexity for maps with values in length spaces} to study the 
convexity properties of the image of the \textbf{cylinder valued momentum map}. This object,
introduced by Condevaux, Dazord, and 
Molino~\cite{condevaux dazord and molino},
naturally generalizes the standard momentum map definition due to Kostant and Souriau. The standard momentum map is associated to a canonical Lie algebra action on a symplectic manifold. Its convexity properties have been extensively studied~\cite{atiyah 82, convexity,
convexity2, hilnebplank, sjamaar}.

Unlike the standard momentum map, the cylinder valued momentum map always exists for any
canonical Lie algebra action. However, the convexity properties of the standard momentum map
cannot be trivially extended to this object because it does not map into a vector space but
into a manifold that is, in general, diffeomorphic to a cylinder. Thus, in order to study the convexity properties of the cylinder valued momentum map the notion of convexity introduced and studied in  
\S\ref{Image convexity for maps with values in length spaces} is necessary. 

\subsection{The cylinder valued momentum map}

In the following paragraphs we quickly review the elementary properties of the cylinder valued momentum map. For more information and for detailed proofs the reader is encouraged to check with~\cite{condevaux dazord and molino} or with Chapter 5 of~\cite{hsr}. 

Let $(M, \omega )$ be a
connected paracompact symplectic manifold and let $\mathfrak{g}$ be a Lie algebra that acts
canonically on $M$. Take the Cartesian product $M \times 
\mathfrak{g}^\ast$  and let
$\pi:M\times\mathfrak{g}^\ast\rightarrow M$ be the projection onto
$M$. Consider
$\pi$ as the bundle map of the trivial principal fiber bundle $(M \times
\mathfrak{g}^\ast, M,
\pi, \mathfrak{g}^\ast)$ that has $(\mathfrak{g}^\ast,+) $ as Abelian structure
group. The group
$(\mathfrak{g}^\ast,+) $ acts on $M \times \mathfrak{g}^\ast $ by $\nu \cdot (m,
\mu):=(m, \mu- \nu)$, with $m \in M $ and $\mu, \nu \in \mathfrak{g}^\ast$. Let
$\alpha
\in \Omega^1(M \times \mathfrak{g}^\ast; \mathfrak{g}^\ast)$ be the connection
one-form defined by
\begin{equation}
\label{definition of alpha connection}
\langle \alpha(m , \mu) (v _m, \nu), \xi\rangle:=
(\mathbf{i}_{\xi_M} \omega) (m)
(v _m) -\langle \nu, \xi \rangle,
\end{equation} where $(m, \mu)\in M \times \mathfrak{g}^\ast $, $(v _m, \nu) \in T
_m M \times
\mathfrak{g}^\ast $,  $\langle\cdot , \cdot \rangle $ denotes the natural pairing
between
$\mathfrak{g}^\ast $ and $\mathfrak{g}$, and $\xi_M $ is the infinitesimal
generator vector field associated to $\xi\in \mathfrak{g}$. The connection 
$\alpha$  is flat.  For $(z,\mu) \in M \times
\mathfrak{g}^\ast$, let $(M
\times
\mathfrak{g}^\ast)(z, \mu) $ be the holonomy bundle through
$(z, \mu) $ and let $ {\mathcal H}(z , \mu)$ be the holonomy group of $\alpha$
with reference point $(z, \mu) $ (which is an Abelian zero dimensional Lie subgroup of
$\mathfrak{g}^\ast$ by the flatness of $\alpha$). The
principal bundle 
$((M
\times
\mathfrak{g}^\ast)(z, \mu),M, \pi|_{(M \times \mathfrak{g}^\ast)(z,
\mu)},{\mathcal H}(z , \mu))
$ is a reduction of the principal bundle $(M \times \mathfrak{g}^\ast, M, \pi,
\mathfrak{g}^\ast)$. To
simplify notation, we will write $(\widetilde{M}, M, \widetilde{p},
{\mathcal H}) $ instead of $((M \times \mathfrak{g}^\ast)(z, \mu),M,
\pi|_{(M \times \mathfrak{g}^\ast)(z, \mu)},{\mathcal H}(z , \mu)) $. Let
$\widetilde{\mathbf{K}}: \widetilde{M} \subset M \times
\mathfrak{g}^\ast\rightarrow \mathfrak{g}^\ast$ be the projection into the
$\mathfrak{g}^\ast$-factor. 

Let $\overline{{\mathcal H}}$ be the closure of ${\mathcal H} $ in
$\mathfrak{g}^\ast$. Since $\overline{{\mathcal H}}$ is a closed subgroup of
$(\mathfrak{g}^\ast, +)$, the quotient $C:= \mathfrak{g}^\ast/ \overline{{\mathcal
H}}$ is a cylinder (that is, it is isomorphic to the Abelian Lie group
$\mathbb{R}^a \times
\mathbb{T}^b$ for some $ a, b \in \mathbb{N}$). Let
$\pi_C:
\mathfrak{g}^\ast\rightarrow
\mathfrak{g}^\ast/\overline{{\mathcal H}}$ be the projection. Define
$\mathbf{K}: M \rightarrow  C $ to be the map that makes the following diagram
commutative:
\begin{equation}
\label{diagram commutative cylinder valued momentum map}
\begin{CD}
\widetilde{M}@>\widetilde{\mathbf{K}}>>\mathfrak{g}^\ast\\ @V\widetilde{p} VV	
@VV\pi_C V\\ M@>\mathbf{K}>>\mathfrak{g}^\ast/ \overline{ {\mathcal H}}.
\end{CD}
\end{equation} In other words, $\mathbf{K}$ is defined by  $ \mathbf{K}(m)= \pi_C
(\nu) $, where
$\nu \in
\mathfrak{g}^\ast $ is any element  such that $(m, \nu) \in \widetilde{M } $.

We will refer to $ \mathbf{K}:M \rightarrow  \mathfrak{g}^\ast/ \overline{{\mathcal
H}}$ as a \textbf{cylinder valued momentum map} associated to the canonical
$\mathfrak{g}$-action on $(M, \omega)$ and to ${\cal H}$ as the 
\textbf{Hamiltonian holonomy}
of the $\mathfrak{g}$-action on $(M , \omega  )$.

\medskip

\noindent \textbf{Elementary properties.} The cylinder valued  momentum  map is a strict
generalization of the standard (Kostant-Souriau) momentum map since  the
$\mathfrak{g}$-action has a standard momentum map if and only if the holonomy group
${\mathcal H} $ is trivial. In such a case the cylinder valued  momentum map  is a
standard momentum map. The cylinder valued momentum map  satisfies Noether's Theorem, that
is, for any $\mathfrak{g}
$-invariant function $h \in  C^\infty(M)^{\mathfrak{g}}:=\{f \in  C^\infty(M)\mid
\mathbf{d}h (\xi _M) = 0\text{ for all }\xi\in \mathfrak{g}\}
$, the flow $F _t 
$ of its associated Hamiltonian vector field $X _h$ satisfies the
identity
$
\mathbf{K} \circ F _t= \mathbf{K}| _{{\rm Dom}(F _t)}$. Additionally,
for any $ v _m \in T _mM $, $m \in M $, $T _m \mathbf{K} ( v _m) = T _\mu \pi_C ( \nu)$,
where $\mu \in \mathfrak{g}^\ast  $ is any element such that $\mathbf{K} (m)=
\pi_C (\mu)
$ and $\nu \in \mathfrak{g}^\ast $ is uniquely determined by
$
\langle \nu, \xi\rangle=(\mathbf{i}_{\xi_M} \omega)(m) (v _m)$,
for any
$\xi\in \mathfrak{g}$.
Also,
$\ker (T _m \mathbf{K})= \left( \left({\rm Lie}(\overline{{\mathcal H}})\right)
^{\circ}\cdot m\right) ^\omega$, where ${\rm Lie}(\overline{{\mathcal H}}) \subset \mathfrak{g}^\ast$ is the Lie algebra of $\overline{{\mathcal H}}$,  and 
${\rm range}\, (T _m \mathbf{K})= T _\mu \pi_C \left((\mathfrak{g}_{m})^\circ
\right)$ (Bifurcation Lemma). In the first statement we use the fact that  
the annihilator $\left({\rm Lie}(\overline{{\mathcal H}})\right)^{\circ}$ is a Lie subalgebra of $\mathfrak{g}$. The notation $\mathfrak{k} \cdot m $ for any Lie subalgebra $\mathfrak{k} \subset \mathfrak{g}$ means the vector subspace of $T _mM $ formed by evaluating all infinitesimal generators $\eta_M$ at the point $m \in M $ for all $\eta \in \mathfrak{k}$. The upper index $\omega$ denotes the $\omega$-orthogonal complement of the set in question.

\medskip

\noindent \textbf{Equivariance properties of the cylinder valued momentum map.}
Suppose now that the $\mathfrak{g}$-Lie algebra action on $(M, \omega)$ 
is obtained from a canonical action
of the Lie group $G$ on $(M, \omega)$ by taking the  infinitesimal
generators of all elements  in $\mathfrak{g}$. There is a
$G$-action on the target space of the cylinder valued momentum map $\mathbf{K}:M
\rightarrow
 \mathfrak{g}^\ast/\overline{{\mathcal H}} $ with respect to which it is
$G$-equivariant. This action is constructed by noticing first that  
the Hamiltonian holonomy
${\cal H}$ is invariant under the coadjoint action, that is,
$
\mbox{\rm Ad}^\ast _{g ^{-1}} {\mathcal H}\subset  {\mathcal H}$, for any $g
\in  G $.
Actually, if $G$ is connected, then ${\mathcal H} $  
is pointwise fixed by the coadjoint action~\cite{cylinder reduction}.
Hence, there is a unique group action $\mathcal{A}d^\ast : G
\times 
\mathfrak{g}^\ast/\overline{{\mathcal H}} \rightarrow  \mathfrak{g}^\ast/
\overline{{\mathcal H}}
$ such that  for any $g \in G $, 
$\mathcal{A}d^\ast_{ g ^{-1}} \circ \pi _C= \pi _C \circ \mbox{\rm Ad} ^\ast  _{g
^{-1}}$. With this in mind, we define
 $\overline{\sigma}: G\times 
M \rightarrow \mathfrak{g}^\ast /\overline{{\mathcal H}}$ by
\begin{equation*}
\label{definition of the non-equivariance cocycle}
\overline{\sigma}(g,m):= \mathbf{K}(\Phi_g( m))-\mathcal{A}d^\ast_{g ^{-1}}
\mathbf{K} (m). 
\end{equation*} 
~Since $M $ is connected by hypothesis, it can be shown that  $ \overline{\sigma} $ does not
depend on the points $m \in  M$ and hence it defines a map  
$\sigma: G
\rightarrow  \mathfrak{g}^\ast /\overline{{\mathcal H}} $ which is a group valued one-cocycle:
for any $g,h
\in G $, it satisfies the equality
$
\sigma(gh)= \sigma(g)+ \mathcal{A}d^\ast _{g ^{-1}} \sigma (h)$. This guarantees that
the map 
\begin{equation*}
\begin{array}{cccc}
\Theta:& G \times  \mathfrak{g}^\ast /\overline{{\mathcal H}}&\longrightarrow &
\mathfrak{g}^\ast /\overline{{\mathcal H}}\\
	&(g, \pi_C(\mu))&\longmapsto&
\mathcal{A}d^\ast _{g ^{-1}}(\pi_C(\mu))+ \sigma (g),
\end{array}
\end{equation*} defines a $G$-action on $\mathfrak{g}^\ast /\overline{{\mathcal
H}}$ with respect to which the cylinder valued  momentum map $\mathbf{K}$ is
$G$-equivariant, that is, for any $g  \in G $, $m \in M $, we have
\[
\mathbf{K}(\Phi_g(m))= \Theta_g(\mathbf{K}(m)).
\]
We will refer to $\sigma: G \rightarrow  \mathfrak{g}^\ast
/\overline{{\mathcal H}} $ as the \textbf{non-equivariance one-cocycle} of the 
cylinder valued    momentum map
$\mathbf{K}:M \rightarrow \mathfrak{g}^\ast/ \overline{{\mathcal H}} $ and to
$\Theta$ as the \textbf{affine $G$-action} on $\mathfrak{g}^\ast
/\overline{{\mathcal H}}$ induced by $\sigma$.
The infinitesimal
generators of the affine $G$-action on $\mathfrak{g}^\ast/\overline{{\mathcal H}}$
are given by the expression
\begin{equation}
\label{infinitesimal generators of the affine}
\xi_{\mathfrak{g}^\ast/ \overline{{\mathcal H}}} (\pi _C (\mu))=-T _\mu \pi _C
\left(\Psi (m)(\xi, \cdot )
\right),
\end{equation} for any $\xi \in  \mathfrak{g}$, $(m, \mu)\in  \widetilde{M} $, 
where $\Psi:M \rightarrow Z ^2(\mathfrak{g})$ is the Chu map defined by 
$\Psi(\xi, \eta):= \omega\left(\xi_M, \eta_M\right)$, for any $\xi,
\eta\in \mathfrak{g}$.

\medskip

\noindent \textbf{The Poisson structure on $\mathfrak{g}^\ast/ \overline{{\mathcal H}}$.}
The bracket
$\{\cdot,\cdot\}_{\mathfrak{g}^\ast/\overline{{\mathcal H}}} :
C^{\infty}(\mathfrak{g}^\ast/\overline{{\mathcal H}})\times
C^{\infty}(\mathfrak{g}^\ast/\overline{{\mathcal H}}) \rightarrow
\mathbb{R} $ defined by 
\begin{equation*}
\label{bracket on target space}
\{ f , g\}_{\mathfrak{g}^\ast/\overline{{\mathcal H}}}(\pi_C (\mu)) =
\Psi (m) \left(\frac{\delta(f \circ \pi _C)}{\delta \mu}, \frac{\delta(g
\circ \pi _C)}{\delta \mu} \right),
\end{equation*} 
where $f , g \in C^{\infty}(\mathfrak{g}^\ast/\overline{{\mathcal H}}) $,
$(m, \mu)\in  \widetilde{M} $, $\pi _C: \mathfrak{g}^\ast \rightarrow
\mathfrak{g}^\ast/\overline{{\mathcal H}} $ is the projection, and $\Psi: M
\rightarrow Z ^2 (\mathfrak{g})$ is the Chu map, defines a Poisson 
structure on $\mathfrak{g}^\ast/ \overline{{\mathcal H}} $ such that 
$\mathbf{K}:M \rightarrow 
\mathfrak{g}^\ast/\overline{{\mathcal H}} $ is a Poisson map.

\subsection{A normal form for the cylinder valued  momentum map}
\label{A normal form for the cylinder valued  momentum map}

A major technical tool in some proofs of the classical convexity theorems  for the standard
momentum map is a normal form  obtained  by Marle~\cite{marle}
and by Guillemin and Sternberg~\cite{GuSt1984}. This normal form is a version of the
classical Slice Theorem for proper group actions adapted to the symplectic symmetric setup
that provides a semi-global set of coordinates (global only in the direction of the group
orbits) in which the standard momentum map takes a particularly convenient and simple form
and in which the conditions of the Local-to-Global
Principle can be verified. This normal form has
been generalized to the context of the cylinder valued momentum map in~\cite{symplectic
slice, stratified cylinder valued momentum map}. We briefly review this generalization in the
following paragraphs.  

In this section we will work on a  connected and paracompact symplectic manifold $(M,
\omega)$ acted properly and symplectically upon by the Lie group $G$ with Lie algebra
$\mathfrak{g}$.  The first step in the construction of the symplectic slice theorem is the
splitting of the Lie algebra $\mathfrak{g}$ of $G$ into three parts. The first summand is
defined by 
\begin{equation}
\label{definition of k subalgebra}
\mathfrak{k}:=\left\{ \xi\in \mathfrak{g}\mid \xi_M (m) \in (\mathfrak{g}\cdot m)^{\omega(m)} \right\},
\end{equation}
where $m \in M$ is the point around whose $G$-orbit we want to construct the symplectic slice.
The set $\mathfrak{k}$ is clearly a vector subspace of $\mathfrak{g} $ that contains
the Lie algebra $\mathfrak{g}_m $ of the isotropy subgroup $G _m $ of the point $m\in M$. In fact, $\mathfrak{k}$ is a Lie subalgebra of 
$\mathfrak{g}$. Since the $G $-action is by hypothesis proper, the isotropy subgroup $G _m $  is compact and hence there is an
$\mbox{\rm Ad} _{G_{m} }$-invariant inner product $\langle \cdot , \cdot \rangle_{\mathfrak{g}} $ on $\mathfrak{g}$. We decompose 
\begin{equation}
\label{splitting of lie algebras for slice}
\mathfrak{k}= \mathfrak{g}_{m}\oplus \mathfrak{m} \quad\text{and}\quad \mathfrak{g}=
\mathfrak{g}_{m}\oplus \mathfrak{m}\oplus \mathfrak{q},
\end{equation}
where $\mathfrak{m}$ is the $\langle \cdot , \cdot \rangle_{\mathfrak{g}} $--orthogonal
complement of
$\mathfrak{g}_{m}$ in $\mathfrak{k}$ and $\mathfrak{q}$ is the $\langle \cdot , \cdot
\rangle_{\mathfrak{g}} $--orthogonal complement of $\mathfrak{k}$ in $\mathfrak{g}$. The
splittings in~(\ref{splitting of lie algebras for slice}) induce similar ones on the duals 
\begin{equation}
\label{splitting of lie algebras for slice dual}
\mathfrak{k}^\ast = \mathfrak{g}_{m}^\ast \oplus \mathfrak{m}^\ast  \quad\text{and}\quad
\mathfrak{g}^\ast =
\mathfrak{g}_{m}^\ast \oplus \mathfrak{m}^\ast \oplus \mathfrak{q}^\ast. 
\end{equation}
Each of the spaces in this decomposition should be understood as the set of covectors in
$\mathfrak{g}^\ast $ that can be written as $\langle \xi, \cdot \rangle _{\mathfrak{g}}$, with
$\xi $ in the corresponding subspace. For example, $ \mathfrak{q}^\ast =\{ \langle \xi, \cdot
\rangle_{\mathfrak{g}}\mid  \xi \in \mathfrak{q}\} $.

The subspace $\mathfrak{q}\cdot m $ is a symplectic subspace of $\left(T _mM, \omega(m)\right)$.

Let now $\ll \cdot , \cdot \gg $ be a $G_{m} $--invariant inner product in $T
_mM $ (available again by the compactness of $G _m  $). Define $V $ as the
orthogonal complement to $\mathfrak{g}\cdot m\cap (\mathfrak{g}\cdot m)^{\omega(m)} =
\mathfrak{k}\cdot m$ in
$(\mathfrak{g}\cdot m) ^{\omega(m)} $ with respect to $\ll \cdot , \cdot \gg $, that is:
\begin{equation*}
(\mathfrak{g}\cdot m) ^{\omega(m)}=\mathfrak{g}\cdot m\cap (\mathfrak{g}\cdot
m)^{\omega(m)} \oplus V= \mathfrak{k}\cdot m \oplus V.
\end{equation*}
The subspace $V  $ is a symplectic $G_{m}  $--invariant subspace of $(T _mM, \omega (m))$
such that $V\cap \mathfrak{q}\cdot m=\{0\} $.
Any such
space $V$ is called a \textbf{symplectic normal space}
at $m$. Since the $G_{m}$--action on $(V, \omega(m)|_{V}) $ is linear and
symplectic it has a standard equivariant
associated momentum map $\mathbf{J}_V:V \rightarrow
\mathfrak{g}_{m} ^\ast$ given by $\left\langle \mathbf{J}_V(v), \eta \right\rangle =
\frac{1}{2} \omega(m)(\xi _V(v), v )$. The proof of the following two results  can be found
in~\cite{symplectic slice, hsr}.

\begin{proposition}[The symplectic tube]
\label{design of symplectic tube at a point}
Let $(M, \omega)$ be a connected paracompact symplectic manifold and $G$  a Lie group acting
properly and canonically on it. Let $m \in M $, $V$ be a symplectic normal space at $m$, and
$\mathfrak{m} \subset \mathfrak{g}$
 the subspace introduced in the splitting~\eqref{splitting of lie algebras for slice}.
Then there exist $G_{m} $--invariant neighborhoods
$\mathfrak{m}^\ast_r$ and $V _r$ of the origin in
$\mathfrak{m}^\ast$ and $V$, respectively, such that the twisted
product
\begin{equation}
\label{marle tube canonical}
Y _r:= G \times_{G_{m} } \left(\mathfrak{m}^\ast_r\times V _r\right)
\end{equation}
is a symplectic manifold  
with the two--form $\omega_{Y _r } $ defined by:
\begin{align}
\omega_{Y _r }&([g, \rho,v])(T_{(g, \rho, v)}\pi(T _e L _g (\xi_1), \alpha_1, u
_1),T_{(g, \rho, v)}\pi(T _e L _g (\xi_2), \alpha_2, u _2))\notag\\
	&:=\langle \alpha_2+ T _v \mathbf{J}_V (u _2), \xi_1\rangle-\langle \alpha_1+ T _v
\mathbf{J}_V ( u _1), \xi_2\rangle+\langle \rho+ \mathbf{J}_V (v),[\xi_1,
\xi_2]\rangle\notag\\
	&\ \ + \Psi (m)(\xi_1, \xi_2)+ \omega(m)(u _1, u _2),
\label{symplectic form in the symplectic tube}
\end{align}
where $\Psi: M
\rightarrow Z ^2 (\mathfrak{g})$ is the Chu map associated to the $G$--action on $(M,
\omega)$,  $\pi:G
\times
\left(\mathfrak{m}^\ast_r\times V _r\right)\rightarrow  G
\times_{G_{m} } \left(\mathfrak{m}^\ast_r\times V _r\right) $ is the projection, $[g,
\rho,v]  \in Y _r$, $\xi_1, \xi_2 \in \mathfrak{g} $, $\alpha_1, \alpha_2 \in
\frak{m}^\ast $, and $u _1, u _2 \in V $.

The Lie group $G$ acts canonically on $(Y _r, \omega_{Y _r})$ by 
$g\cdot [h, \eta, v]:=[gh, \eta,v]$, for any $g \in G $ and any $[h,
\eta, v] \in Y _r$. 
\end{proposition}

In the sequel will refer to the symplectic manifold $(Y _r, \omega_{Y _r})$ as a \textbf{symplectic tube} 
of $(M, \omega)$ at the point $m$.

\begin{theorem}[Symplectic Slice Theorem]
\label{normal form for canonical actions}
Let $(M,\omega)$ be a symplectic manifold and let $G$ be a Lie
group acting properly and canonically on $M$. Let $m\in M$ and let $(Y _r, \omega_{Y
_r}) $ be the $G$--symplectic tube at that point constructed in Proposition~{\rm
\ref{design of symplectic tube at a point}}.  Then there is a $G$--invariant neighborhood
$U$ of $m$ in $M$ and a $G$--equivariant symplectomorphism $\phi:U\rightarrow Y_r$
satisfying
$\phi(m)=[e,\,0,\,0]$.
\end{theorem}

We now provide an expression in the symplectic tube  for the cylinder valued
momentum map. This is what we call the normal form for the cylinder valued   momentum map. The proof of the
following theorem can be found in~\cite{stratified cylinder valued momentum map}. 

\begin{theorem}[Normal form for the cylinder valued  momentum map]
\label{normal form theorem for cylinder valued momentum map}
Let $(M, \omega)$ be a connected paracompact symplectic manifold acted properly and
canonically upon by the connected Lie group $G$. Let $m \in M $ and $(Y
_r, \omega _{Y _r})$ be a symplectic tube at $m$ that models a $G$--invariant neighborhood
$U$of the orbit $G \cdot m  $ via the $G$--equivariant symplectomorphism $\phi:(Y
_r, \omega _{Y _r}) \rightarrow (U, \omega| _U)$. Let $\mathbf{K}:M \rightarrow
\mathfrak{g}^\ast/ \overline{{\mathcal H}} $ be a cylinder valued momentum map
associated to the $G$--action on $M$ with non-equivariance one-cocycle $\sigma: G \rightarrow
\mathfrak{g}^\ast/ \overline{{\mathcal H}} $. Then for any $[g, \rho, v] \in
Y _r$ we have 
\begin{eqnarray}
\mathbf{K}(\phi[g, \rho, v])&=& \Theta_g\left(\mathbf{K}(m)+\pi_C(\rho+
\mathbf{J}_V (v)\right)\label{normal form cylinder valued
momentum map}\\
	&= & \Theta_g\left(\mathbf{K}(m)\right)+\pi_C\left(\mbox{\rm Ad}^\ast _{g ^{-1}}(\rho+
\mathbf{J}_V (v))\right)\label{normal form cylinder valued
momentum map 1}
\end{eqnarray} 
where $\pi_C: \mathfrak{g}^\ast  \rightarrow \mathfrak{g}^\ast/ \overline{ {\mathcal H}} $
is the projection and $\Theta:G \times  \mathfrak{g}^\ast/ \overline{{\mathcal H}} \rightarrow 
\mathfrak{g}^\ast/ \overline{{\mathcal H}} $ is the affine action associated to the
non-equivariance one-cocycle $\sigma$.
\end{theorem}

\subsection{Closed Hamiltonian holonomies and covering spaces}

We will use in our study of the convexity properties of the cylinder valued momentum map a
hypothesis that allows us to naturally endow the target space of this map with all the
necessary metric properties . More specifically, \emph{we will assume in all that follows
that the Hamiltonian holonomy ${\mathcal H} $ is a closed subgroup of $(\mathfrak{g}^\ast, +)$.}

In order to spell out the implications of this hypothesis we introduce some terminology. Let 
$G$ be a group that acts on a topological space $X$. This action is
called \textbf{totally discontinuous} if every point $x \in X$ has a neighborhood
$U$ such that $g \cdot U\cap U=\varnothing$ for all $g\in G$ satisfying
$g\cdot x\neq x$. 

Let $X$ and $Y$ be two topological spaces and $f:X\rightarrow Y$ a
continuous map. An open set $V\subset Y$ is said to be \textbf{evenly
covered} if it its inverse image $f^{-1}(V)$ is a disjoint union of
open sets $U_i\subset X$ such that the restrictions $\left. f\right|_{
U_i}:U_i\rightarrow V$ are homeomorphisms. The map $f$ is a
\textbf{covering map} if every point $y\in Y$ has an evenly covered
neighborhood.

A continuous map $f:X\rightarrow Y$ between two topological spaces
induces a homomorphism of fundamental groups
$f_{*}:\pi_{1}(X,x_0)\rightarrow \pi_{1}(Y,y_0)$ which maps the
class of a loop $\gamma$ in $X$ to the class of the loop
$f\circ\gamma$ in $Y$. A covering map $f:X\rightarrow Y$ is said
to be \textbf{regular} if the following equivalent properties hold:

\begin{description}
\item [(i)]  $f_{*}(\pi_{1}(X,x_{0}))$ is a normal subgroup of
$\pi_{1}(Y,y_{0})$.
\item [(ii)]  $f_{*}(\pi_{1}(X,x_{0}))$ does not depend on $x_{0}\in
f^{-1}(x_{0})$.
\end{description}
A \textbf{deck transformation} of a covering $f:X\rightarrow Y$ is a
homeomorphism which permutes the sheets of the cover or
equivalently permutes the points of $f^{-1}(y)$ for any $y\in Y$. The proof of the
following two results can be found in~\cite{burago}, Propositions 3.4.15 and 3.4.16. 

\begin{proposition}  
\label{g invariance projection}
Let $G$ be a group acting on a
topological space $X$ freely and totally discontinuously. Then, the
projection $\pi_{G}:X\rightarrow X/G$ onto the orbit space is a regular covering map.
Moreover,  the group of its deck
transformations coincides with $G$.
\end{proposition}

\begin{theorem} 
\label{theorem for deck}
Let $f:X\rightarrow Y$ be a regular
covering and $G$ its group of deck transformations. Then, the
length metrics on $Y$ are in one-to-one correspondence with the
$G$-invariant length metrics on $X$ so that for corresponding
metrics $d_X$ on $X$ and $d_Y$ on $Y$, $f$ is a local isometry.
\end{theorem}

\begin{lemma}
\label{lemma for completeness}
Let $X$ and $Y$ be two metric spaces and $f:X\rightarrow Y$ a
covering map that is also a local isometry. If $X$ is complete
then $Y$ is complete.
\end{lemma}

\noindent\textbf{Proof.\ \ } Let $\{y_n\}_{n\in \mathbb{N}}$ be a Cauchy sequence in
$(Y,d_Y)$. Given that the map $f$ is  a covering map and a local isometry there 
exists $\varepsilon_0 > 0$  small enough so that
 $f^{-1}(B_Y(y,\varepsilon_0))=\bigcup_{x\in
f^{-1}(y)}B_X(x,\varepsilon_0)$ is a disjoint union. As $\{y_n\}_{n\in \mathbb{N}}$ is a
Cauchy sequence, there exists
$N_{\varepsilon_0}\in \mathbb{N}$ such that
$d_Y(y_n,y_m)<N_{\varepsilon_0}$ for all $n,m\geq
N_{\varepsilon_0}$. Consequently, for $m\geq N_{\varepsilon_0}$ we
have that $y_m\in B_Y(y_{_{N_{\varepsilon_0}}},\varepsilon_0)$. 
Select an
arbitrary ball $B_X(x_0,\varepsilon_0)$ in
the disjoint union $\bigcup_{x\in
f^{-1}(y_{_{N_{\varepsilon_0}}})}B_X(x,\varepsilon_0)$. 
Define the sequence $\{x_m\}_{m\geq N_{\varepsilon_0}}$, $x_m:=f^{-1}(y_m)\cap
B_X(x_0,\varepsilon_0)$. Since $f$ is a local isometry, the sequence $\{x_m\}_{m\geq N_{\varepsilon_0}}$ is clearly  Cauchy
 in
$X$ and hence the completeness of $X$ and the continuity
of $f$ ensures that $\{y_n\}_{n\in \mathbb{N}}$ has a convergent
subsequence.
\quad $\blacksquare$

\medskip

We now come back to the implications of the closedness hypothesis on the Hamiltonian holonomy
${\mathcal H}$.

\begin{proposition}
\label{properties of pic closed h}
Let $(M, \omega)$ be a connected paracompact symplectic manifold acted 
canonically upon by the Lie algebra $\mathfrak{g}$ with Hamiltonian holonomy ${\mathcal H}
\subset \mathfrak{g}^\ast$. If ${\mathcal H} $ is a closed subset of $\mathfrak{g}^\ast$
then: 
\begin{description}
\item [(i)] The projection $\pi _C: \mathfrak{g}^\ast \rightarrow  \mathfrak{g}^\ast/
{\mathcal H}
$ is a regular covering smooth map and hence the Euclidean metric in $\mathfrak{g}^\ast$
projects naturally to a length metric on $\mathfrak{g}^\ast/ {\mathcal H}  $ with respect
to which this space is complete and locally compact and  the projection $\pi _C $
is  a local isometry.
\item [(ii)] Suppose that there exists a compact connected Lie group $G$ whose Lie
algebra is
$\mathfrak{g} $. Identify the positive Weyl chamber $\mathfrak{t}^\ast _+ \subset
\mathfrak{g}^\ast$ with the orbit space of the coadjoint action of $G$ on
$\mathfrak{g}^\ast$. Then, $\mathfrak{t}^\ast _+ $ is a closed convex subset of
$\mathfrak{g}^\ast $ and hence it has a natural length space structure with respect to which
it is complete and locally compact.
${\mathcal H}$ acts on
$\mathfrak{t}^\ast _+$ in a totally discontinuous fashion and hence the length metric in
$\mathfrak{t}^\ast _+$ projects naturally to a length metric on $\mathfrak{t}^\ast_+/
{\mathcal H}  $ with respect to which this space is complete and locally compact and  the
orbit space projection $\pi _C^+ :\mathfrak{t}^\ast _+ \rightarrow  \mathfrak{t}^\ast _+/
{\mathcal H}$ is  a local isometry.
\item [(iii)]  In the hypotheses of part \textbf{(ii)}, the natural identification of
$\mathfrak{t}^\ast _+ $ with the orbit space $\mathfrak{g}^\ast /G $ of the coadjoint action
induces an identification of $\mathfrak{t}^\ast _+/ {\mathcal H} $ with the orbit space
$(\mathfrak{g}^\ast/{\mathcal H})/G $ of the $\mathcal{A}d^\ast $ action of $G$ on 
$\mathfrak{g}^\ast/ {\mathcal H}$ with respect to which the following diagram commutes: 
\begin{equation}
\label{commuting pi}
\begin{CD}
\mathfrak{g}^\ast@>\pi _C>>\mathfrak{g}^\ast/ {\mathcal H}\\
@V\pi_G VV		@VV\pi_G^+ V\\
\mathfrak{g}^\ast/G\simeq \mathfrak{t}^\ast _+@>\pi_C^+>>\mathfrak{t}^\ast
_+/ {\mathcal H}\simeq (\mathfrak{g}^\ast/  {\mathcal H})/G.
\end{CD}
\end{equation}
\end{description}
\end{proposition} 
 
\noindent\textbf{Proof.\ \ } Since ${\mathcal H}$ acts on $(\mathfrak{g}^\ast,+)$ by
translations, the Euclidean metric in $\mathfrak{g}^\ast $ is ${\mathcal H} $-invariant.
Additionally, as this action is free and proper, the Slice Theorem guarantees that any point
$\mu\in \mathfrak{g}^\ast$ has a ${\mathcal H}$-invariant neighborhood that is equivariantly
diffeomorphic to the product ${\mathcal H} \times U $, with $U \subset \mathfrak{g}^\ast$ an
open neighborhood of zero in $\mathfrak{g}^\ast$. In this semiglobal model the point
$\mu$ is represented by the element $(0,0)$. Since ${\mathcal H}$ is a closed zero
dimensional submanifold of $\mathfrak{g}^\ast$, it follows that the set
$\{0\}\times  U $ is a open neighborhood of $(0,0)\equiv \mu$. Moreover, for any $\nu \in
{\mathcal H} $ different from zero, we have  $\nu \cdot (\{0\}\times  U)=\{\nu\}\times  U  $
and since $\left(\{\nu\}\times  U\right) \cap \left(\{0\}\times  U \right)=\varnothing $ we conclude that ${\mathcal
H} $ acts totally discontinuously on $\mathfrak{g}^\ast$. The statement in part \textbf{(i)}
follows then by Proposition~\ref{g invariance projection}, Theorem~\ref{theorem for deck},
and Lemma~\ref{lemma for completeness}. The local compactness of
$\mathfrak{g}^\ast/{\mathcal H}$ is a consequence of the open character of the orbit
projection $\pi _C $ and the local compactness of $\mathfrak{g}^\ast$ (we recall that
every orbit projection is an open map and that  if $f:X\rightarrow Y$ is an arbitrary open
map from a locally compact topological space onto a Hausdorff space $Y$, then
$Y$ is locally compact).

As to part \textbf{(ii)} we recall that the positive Weyl chamber $\mathfrak{t}^\ast _+ $
is  a closed convex subset of $\mathfrak{g}^\ast$ and hence a complete and locally
compact length metric space (see Definition~\ref{metric convex} and
Lemma~\ref{caracterizare}). We now recall that if $G$ is connected, then ${\mathcal H} $  
is pointwise fixed by the coadjoint action~\cite{cylinder reduction} and hence the
$G$-coadjoint action  and the  ${\mathcal H}$-action on $\mathfrak{g}^\ast$ commute, which
guarantees that the ${\mathcal H}$-action on $\mathfrak{g}^\ast$ drops to an ${\mathcal
H}$-action on $\mathfrak{g}^\ast/G\simeq \mathfrak{t}^\ast _+ $. Since this action can be
viewed as the restriction to $\mathfrak{t}^\ast _+ $ of the totally discontinuous ${\mathcal
H}$-action on $\mathfrak{g}^\ast$, we conclude that it is also totally discontinuous and
hence Proposition~\ref{g invariance projection}, Theorem~\ref{theorem for deck}, and
Lemma~\ref{lemma for completeness} apply, which establishes the statement.

Regarding part \textbf{(iii)}, the identification $\mathfrak{t}^\ast
_+/ {\mathcal H}\simeq (\mathfrak{g}^\ast/  {\mathcal H})/G $ is a consequence of the fact
that the
$G$-coadjoint action  and the  ${\mathcal H}$-action on $\mathfrak{g}^\ast$ commute because
$G$ is connected. The rest of the statement is a straightforward diagram chasing exercise.
\quad $\blacksquare$

\subsection{Convexity properties of the cylinder valued momentum map. The Abelian case.}

In this subsection it will be shown that the image of the cylinder valued momentum map associated
to a proper Abelian Lie group action is weakly convex. The approach taken to
prove this statement consists in using the Symplectic Slice Theorem~\ref{normal form for canonical 
actions} to show that this map satisfies the local hypotheses needed to apply  the generalization of the
Local-to-Global Principle for length spaces (Theorem~\ref{lokal to global metric}). The main step
in that direction is taken in the following proposition.

\begin{proposition}
\label{local expression j u}
Let $(M,\omega)$ be a connected paracompact symplectic manifold and let $G$ be a connected
Abelian Lie group acting properly and canonically on $M$ with closed Hamiltonian holonomy
${\mathcal H}$. Let $\mathbf{K}:M \rightarrow \mathfrak{g}^\ast/ {\mathcal H}$ be a cylinder
valued momentum map for this action and $m  \in M $ arbitrary such that $\mathbf{K}(m)=[
\mu] \in  \mathfrak{g}^\ast/ {\mathcal H} $. Then there exists an open neighborhood $U$ of 
$m$ in $M$ and  open neighborhoods $W $ and $V$ of $\mu\in \mathfrak{g}^\ast$ and  $[\mu]\in
\mathfrak{g}^\ast / {\mathcal H}$, respectively, such that  $\mathbf{K} (U) \subset  V $,
$\left.\pi _C\right|_{W}:W
\rightarrow V $ is a diffeomorphism, and 
\begin{equation}
\label{local expression momentum}
\left.\pi _C\right|_{W}^{-1}\circ  \mathbf{K}|_{U}= \mathbf{J} _U+c,
\end{equation}
with $c \in \mathfrak{g}^\ast$ a constant and  $\mathbf{J}_U:U \rightarrow 
\mathfrak{g}^\ast$ a map that in symplectic slice coordinates around the point $m$ has the
expression
\begin{equation}
\label{local expression J}
\mathbf{J}_U([g, \rho, v])=\rho+ \mathbf{J}_V (v)-\langle
\mathbb{P}_{\mathfrak{q}}(\exp  ^{-1}(s([g]))), \cdot \rangle_{\mathfrak{q}}.
\end{equation}
The neighborhood $U$ has been chosen so that it
can be written in slice coordinates as $U\equiv U _e \times _{G _m}(\mathfrak{m}^\ast \times
 V _r)$, with $U _e $ an open $G _m$-invariant neighborhood of  $e$ in $ G $ small enough so
that there exists a local section $s:U _e/G _m \rightarrow V _e $ for the projection $G
\rightarrow G/G _m$.  $V _e $ is an open neighborhood of  $e \in G $ such that $\exp :U _0
\rightarrow V _e$ is a diffeomorphism, for some open neighborhood $U _0$  of $0 \in
\mathfrak{g}^\ast$.
$\mathbb{P}_{\mathfrak{q}}: \mathfrak{g}=
\mathfrak{g}_{m}\oplus \mathfrak{m}\oplus \mathfrak{q} \rightarrow \mathfrak{q}$ is the
projection onto $\mathfrak{q}$ constructed using the splitting in 
\eqref{splitting of lie algebras for slice} and $\langle\cdot , \cdot \rangle_{\mathfrak{q}}$ is the non-
degenerate bilinear form
on $\mathfrak{q}$ induced by the Chu map at $m$, that is, for any $\xi,
\eta \in  \mathfrak{q}$, $\langle\xi , \eta \rangle_{\mathfrak{q}}:= \omega (m)(\xi_M(m),
\eta_M(m))$.
\end{proposition}  

\noindent\textbf{Proof.\ \ } Due to the closedness hypothesis on the Hamiltonian holonomy
${\mathcal H} $, the projection
$\pi_C$ is a local diffeomorphism (see Proposition~\ref{properties of pic closed h}) and
hence there exists an open neighborhood $V$ of  $\mathbf{K} (m)$ in $\mathfrak{g}^\ast/
{\mathcal H} $ and a neighborhood $W$ of some element in the fiber $\pi _C
^{-1}(\mathbf{K} (m) ) \subset \mathfrak{g}^\ast$ such that
$\left.\pi _C\right|_{W}:W
\rightarrow V $ is a diffeomorphism. Let $U$ be the connected component containing $m$ of
the intersection of $\mathbf{K}^{-1}(V) $ with the domain of a symplectic slice chart around
$m$ and shrink it, if necessary, so that the group factor in the slice coordinates has the
properties in the statement of the proposition. 

We start by noticing that~(\ref{local expression J}) is well defined because for any $h \in
G _m $ and any $[g, \rho, v]\in  U _e \times _{G _m}(\mathfrak{m}^\ast \times
 V _r)$ 
\begin{eqnarray*}
\mathbf{J}_U([gh, h ^{-1}\cdot \rho,h ^{-1}\cdot  v])&=&c+ \mbox{\rm Ad}^\ast _h\rho+
\mathbf{J}_V (h ^{-1}\cdot v)-\langle
\mathbb{P}_{\mathfrak{q}}(\exp  ^{-1}(s([gh]))), \cdot \rangle_{\mathfrak{q}}\\
	&= &c+ \mbox{\rm Ad}^\ast _h(\rho+
\mathbf{J}_V ( v))-\langle
\mathbb{P}_{\mathfrak{q}}(\exp  ^{-1}(s([g]))), \cdot \rangle_{\mathfrak{q}}\\
	&= &c+ \rho+ \mathbf{J}_V (v)-\langle
\mathbb{P}_{\mathfrak{q}}(\exp  ^{-1}(s([g]))), \cdot \rangle_{\mathfrak{q}}.
\end{eqnarray*}
The last equality follows from the Abelian character of $G$.

Next, we will show that for any $z \in U $ and any $\xi\in \mathfrak{g} $, the map
$\mathbf{J}_U ^\xi:=\langle \mathbf{J} _U, \xi\rangle $ satisfies 
\begin{equation}
\label{generator local momentum}
X_{\mathbf{J}_U ^\xi}(z)= \xi_M (z),
\end{equation}
with $X_{\mathbf{J}_U ^\xi} $ the Hamiltonian vector field associated to the function
$\mathbf{J}_U ^\xi \in  C^\infty(M) $. Since this is a local statement, it suffices to show
that
\begin{equation}
\label{thing to verify momentum}
\mathbf{i}_{\xi_{Y _r}} \omega _{Y _r}([\exp \zeta, \rho, v])= \mathbf{d}
\mathbf{J}_{U}^\xi([\exp \zeta, \rho, v]),
\end{equation}
where $[\exp \zeta, \rho, v]$ is the expression of $z$ in slice coordinates and $\zeta \in
\mathfrak{g}$ chosen so that $s([\exp \zeta])=\exp \zeta $. We prove~(\ref{thing to verify 
momentum}) by using the expression of $\omega_{Y _r}$ in Proposition~\ref{design of symplectic tube 
at a point}. First of all, notice that 
\[
\xi_{Y _r}([\exp \zeta, \rho, v])= T_{(\exp \zeta,
\rho, v)} \pi(T _e L _{\exp \zeta}( \xi),0,0),
\]
where $\pi: G \times  (\frak{m}^\ast_r \times V _r)\rightarrow G \times _{G_{m}}(\frak{m}^\ast_r \times 
V _r) $ is the orbit projection. If 
$w := T_{(\exp \zeta,
\rho, v)} \pi(T _e L _{\exp \zeta}( \eta),\alpha,u) \in T_{[\exp \zeta, \rho,v]}(G \times _{G_{m}
}(\frak{m}^\ast_r \times V _r))$ we have
\begin{eqnarray}
\omega_{Y _r}([\exp \zeta, \rho, v])(\xi_{Y _r}([\exp \zeta, \rho, v]),w) &= &\langle
\alpha+T _v \mathbf{J}_V (u), \xi\rangle+ \omega (m)(\xi_M (m), \eta_M (m))\notag\\
	&= &\langle
\alpha+T _v \mathbf{J}_V (u), \xi\rangle+ \omega (m)((\mathbb{P}_{\mathfrak{q}}\xi)_M (m),
( \mathbb{P}_{\mathfrak{q}}\eta)_M (m))\notag\\
	&= &\langle
\alpha+T _v \mathbf{J}_V (u), \xi\rangle+ \langle\mathbb{P}_{\mathfrak{q}}\xi,
\mathbb{P}_{\mathfrak{q}}\eta\rangle_ \mathfrak{q}.\label{first equality symplectic form}
\end{eqnarray}
On the other hand, by \eqref{local expression J}, we have
\begin{equation}
\label{first equality symplectic momentum 2}
\mathbf{d}\mathbf{J}_{U}^\xi([\exp \zeta, \rho, v]) \cdot  w =\langle\alpha+T _v
\mathbf{J}_V (u), \xi\rangle- \left.\frac{d}{dt}\right|_{t=0}\langle \mathbb{P}_
\mathfrak{q} \exp ^{-1} s([\exp\zeta\exp t \eta]), \mathbb{P}_ \mathfrak{q} \xi\rangle_
\mathfrak{q}.
\end{equation}
In order to compute the second summand of the right hand  side, notice that
\[
s([\exp\zeta\exp t \eta])=s([\exp(\zeta+ t \eta)])=\exp(\zeta+ t \eta)h (t),
\]
with $h (t)  $ a curve in $G _m $ such that  $h (0)= e $ and $h' (0)= \lambda \in 
\mathfrak{g}_m $. Consequently,
\[
\left.\frac{d}{dt}\right|_{t=0}s([\exp\zeta\exp t \eta])=T _e L_{\exp \zeta}(\eta+
\lambda)= \left.\frac{d}{dt}\right|_{t=0}\exp (\zeta+t(\eta+ \lambda))
\]
and hence, since $\mathbb{P}_ \mathfrak{q} \lambda =0 $, we have 
\begin{equation}
\label{calculation with p q}
\left.\frac{d}{dt}\right|_{t=0}\langle \mathbb{P}_
\mathfrak{q} \exp ^{-1} s([\exp\zeta\exp t \eta]), \mathbb{P}_ \mathfrak{q} \xi\rangle_
\mathfrak{q}= \left.\frac{d}{dt}\right|_{t=0}\langle \mathbb{P}_
\mathfrak{q} (\zeta+t( \eta+ \lambda)), \mathbb{P}_ \mathfrak{q} \xi\rangle_
\mathfrak{q}=\langle \mathbb{P}_ \mathfrak{q} \eta,\mathbb{P}_ \mathfrak{q} \xi\rangle_
\mathfrak{q}.
\end{equation}
The equalities~(\ref{first equality symplectic form}),~(\ref{first equality symplectic
momentum 2}), and~(\ref{calculation with p q})   show that~(\ref{generator local momentum}) holds. 

With this in mind we will now show that for any $z \in  U $
\begin{equation}
\label{local expression momentum infinitesimal}
T _z(\left.\pi _C\right|_{W}^{-1}\circ  \mathbf{K}|_{U})= T _z\mathbf{J} _U.
\end{equation}
Indeed, for any $v _z\in T _zM $ and $\rho \in \pi _C ^{-1}(V) $ such that  $\mathbf{K}(z)= \pi_C (\rho) 
$,
\[
T _z(\left.\pi _C\right|_{W}^{-1}\circ  \mathbf{K}|_{U})( v_z)
 =\left(T_{\mathbf{K}(z)}\left.\pi _C\right|_{W}^{-1} \circ T _z
\mathbf{K}\right) (v _z) 
= \left(T_{\mathbf{K}(z)}\left.\pi _C\right|_{W}^{-1}\circ  T
_\rho\pi_C\right)(\nu)= \nu,
\] 
where $\nu\in \mathfrak{g}^\ast$ is uniquely determined by the expression
\begin{equation}
\label{determines nu}
\langle \nu, \xi\rangle=(\mathbf{i}_{\xi_{M}} \omega) (z) \left(v _z \right), \text{ for all }
\xi\in  \mathfrak{g}.
\end{equation}
On the other hand, by~(\ref{generator local momentum}), we can write
\[
\langle T _z \mathbf{J} \left(v _z \right), \xi\rangle= \mathbf{d} \mathbf{J}_U ^\xi(z) \left(v_z \right)
=(\mathbf{i}_{\xi_{M}} \omega) (z)\left(v _z\right).
\] 
This, together with~(\ref{determines nu}),  shows that~(\ref{local expression momentum
infinitesimal})  holds. 

Let $c (t) $ be a smooth curve such that $c (0)=m $ and $c (1)=z $,
available by the connectedness of $U$. Then by~(\ref{local expression momentum infinitesimal}) 
\begin{eqnarray*}
\left(\left.\pi _C\right|_{W}^{-1} \circ \mathbf{K} \right)(z) &- &\left(\left.\pi
_C\right|_{W}^{-1} \circ \mathbf{K}\right)(m)
=\int_0^1 \frac{d}{dt} \left(\left.\pi
_C\right|_{W}^{-1}\circ  \mathbf{K}\right)(c (t)) d t\\
	&= & \int_0^1 T_{c (t)}  \left(\left.\pi
_C\right|_{W}^{-1}\circ  \mathbf{K}\right) \left(\dot{c} (t)\right) d t
=\int_0^1 T_{c(t)} \mathbf{J} _U\left(\dot{c} (t)\right) d t
= \mathbf{J} _U (z)- \mathbf{J} _U (m).
\end{eqnarray*}
Since $z \in M $ is arbitrary and $m \in  M $ is fixed, the previous equality shows
that~(\ref{local expression momentum}) holds by setting $c= \left(\left.\pi
_C\right|_{W}^{-1} \circ \mathbf{K}\right)(m)-\mathbf{J} _U (m)$.  
\quad $\blacksquare$

\begin{theorem}
Let $(M,\omega)$ be a connected paracompact symplectic manifold and let $G$ a connected
Abelian Lie group acting properly and canonically on $M$ with closed Hamiltonian holonomy
${\mathcal H}$. Let $\mathbf{K}:M \rightarrow \mathfrak{g}^\ast/ {\mathcal H}$ be a cylinder
valued momentum map for this action. If $\mathbf{K}$ is a closed map then the image
$\mathbf{K}(M) \subset \mathfrak{g}^\ast/{\mathcal H} $  is a weakly convex subset of
$\mathfrak{g}^\ast/ {\mathcal H} $. We think of $\mathfrak{g}^\ast/{\mathcal H}$ as a length
metric space with the length metric naturally inherited from $\mathfrak{g}^\ast$ {\rm (}see
Proposition~{\rm \ref{properties of pic closed h})}. If, in addition, $\mathfrak{g}^\ast/
{\mathcal H} $ is uniquely geodesic then $\mathbf{K} (M) $ is convex, $\mathbf{K} $ has
connected fibers, and it is open onto its image. 
\end{theorem}

\noindent\textbf{Proof.\ \ }We will establish this result by using the
Local-to-Global
Principle for length spaces (Theorem~\ref{lokal to global metric}). First of all,
notice that the closedness of the Hamiltonian holonomy implies, by
Proposition~\ref{properties of pic closed h},  that $\mathfrak{g}^\ast/{\mathcal H}$ is a
complete and locally compact length space and that the projection $\pi_C$ is a local
isometry. Therefore, in order to apply  Theorem~\ref{lokal to global metric} we need to show
that $\mathbf{K} $ is locally open onto its image, locally fiber connected, and has local
convexity data. Now, by Proposition~\ref{local expression j u} (more specifically
by~(\ref{local expression momentum})) and Lemma~\ref{convex image}, it suffices to prove
that those three local properties are satisfied by the map
$\mathbf{J} _U: U \rightarrow  \mathfrak{g}^\ast$.

We start the proof of this fact by recalling that since the $G$-action is proper the
isotropy subgroup $G _m $ is compact and hence its connected component containing the
identity is isomorphic to a torus. Consequently, the  map $\mathbf{J}_V:V
\rightarrow \mathfrak{g}_m^\ast $ is the momentum map of the symplectic representation of a
torus on the symplectic vector space $V $ and hence it automatically has (see for
instance~\cite{hilnebplank} for a proof) local convexity data and it is locally fiber
connected and locally open onto its image. Additionally, if we split $\mathfrak{g}^\ast$ as
$\mathfrak{g}^\ast =\mathfrak{g}_{m}^\ast \oplus \mathfrak{m}^\ast \oplus \mathfrak{q}^\ast$
then the map $\mathbf{J}_U-c $ can be decomposed as
\[
\mathbf{J}_U([g, \rho,v])-c=(\mathbf{J}_V (v), \rho, -\langle
\mathbb{P}_{\mathfrak{q}}(\exp  ^{-1}(s([g]))), \cdot \rangle_{\mathfrak{q}}).
\]
Each of the three components of the map has local convexity data, is locally
open onto its image, and is locally fiber connected. Thus $\mathbf{J}_U$ also has these properties.
\quad $\blacksquare$

\subsection{Convexity for Abelian Lie group valued momentum maps}
\label{Convexity for Abelian Lie group valued momentum maps}

We now discuss the convexity properties of the {\bf  Lie group valued momentum maps} 
introduced in~\cite{mcduff 1988, Ginzburg1992, Huebschmann and Jeffrey 1994, huebschmann
1995, AMM}. We give the definition of  these objects only for
Abelian symmetry groups because in the non-Abelian case these momentum
maps are defined on spaces that are not symplectic (they are
called {\bf  quasi Hamiltonian spaces}) 
thereby leaving the category on which we focus in this paper.

\begin{definition}
Let $(M,\omega)$ be a symplectic manifold and $T ^k$ a torus  acting canonically 
on $(M, \omega)$. The map $\mu:M \rightarrow T^k$ is called a Lie group
valued momentum map if for any $\xi \in \mathfrak{t}$, \[\mathbf{i}_{\xi_M}\omega=\mu^{*}(\theta,
\xi),\] 
where
$\theta\in \Omega^1(T^k,\mathfrak{t})$ is the bi-invariant
Maurer-Cartan form.
\end{definition}

A typical example for this momentum map is provided by the following situation. Take the symplectic 
manifold 
$T^2=S^1\times S^1$ with symplectic form  the standard area form and  consider the action of the 
circle
 on the first circle of $T ^2$. The $S ^1$-valued momentum map associated to this action is the
projection on the second circle of $T ^2$, namely, 
$\mu(e^{i\varphi_1},e^{i\varphi_2})=e^{i\varphi_2}$.

\begin{proposition}
The image of any Abelian Lie group valued momentum map $\mu: M \rightarrow T ^k $ is a  weakly 
convex subset of $T ^k $.
\end{proposition}

\noindent\textbf{Proof.\ \ }Alekseev, Malkin, and
Meinrenken \cite{AMM} (see Proposition 3.4 and Remark 3.3 of this paper) prove that for each point  in 
$M $ there is an open simply connected neighborhood $U \subset M $ and a standard momentum map 
$\Phi:U\rightarrow \mathfrak{t}$ ($\mathfrak{t}$ and $\mathfrak{t}^\ast$ are identified here) such that 
the restriction 
$\exp |_{\Phi (U)}:\Phi (U) \rightarrow T ^k$ is a diffeomorphism onto its image and one has $\mu|_U=
\exp \circ {\Phi}$. This immediately implies that $\mu$ has the local properties (LFC), (LOI), and (LCD). 
Thus the hypotheses of Theorem \ref{lokal to global metric} hold and the statement is a direct 
consequence of this theorem.
\quad $\blacksquare$
\medskip

For a map with values in a metric spaces that is not uniquely geodesic, 
the convexity property of its image is not related to the connectedness of its fibers. This is in sharp 
contrast to the situation encountered for standard momentum maps. For instance, if in the above 
example we multiply the symplectic form by two then the associated $S ^1$-valued
momentum map is $\mu(e^{i\varphi_1},e^{i\varphi_2})=e^{2i\varphi_2}$ which satisfies the hypothesis 
of Theorem~\ref{lokal to global metric} and hence has a  convex image but does not have connected 
fibers.

\subsection{Convexity properties of the cylinder valued momentum map. The Non-Abelian case.}

The study of the convexity properties of the image of the cylinder valued momentum map for
non-Abelian groups presents two main complications with respect to its Abelian analog. First,
unless we assume additional hypotheses, we do not have available a result
similar to Proposition~\ref{local expression j u} that provides a convenient local
representation for $\mathbf{K}$ out of which one can easily conclude the necessary local
properties to ensure convexity out of the Local-to-Global
Principle. Second, the entire image is
not likely to be convex since, already in the standard momentum map case one has to take a
convex piece of the dual of the Lie algebra to obtain convexity. We will take
care of the first problem by working  with special actions, namely those that are
\textbf{tubewise Hamiltonian}, whose definition will be recalled below. As to the second question we will look, not at the image of the cylinder valued momentum map but, as expected from the classical case, at
the intersection of this image with
$\mathfrak{t}^\ast _+ /{\mathcal H}$, taking advantage at the same time of the good behavior of the
projection $\pi _C^+ :\mathfrak{t}^\ast _+ \rightarrow  \mathfrak{t}^\ast _+/
{\mathcal H}$ introduced and discussed in the second part of Proposition~\ref{properties of pic closed
h}.  

\begin{definition}
Let $(M, \omega)$ be a symplectic manifold acted canonically upon by
a Lie group
$G$. For any point $m \in M $, we say that the $G$-action on $M$ is
\textbf{tubewise Hamiltonian at $m$} if there exists a $G$-invariant
open neighborhood of the orbit $G \cdot m $ such that the restriction
of the action to the symplectic manifold $(U, \omega| _U) $ 
has an associated standard momentum map. The $G$-action on $M$ is called \textbf{tubewise Hamiltonian}
if it is tubewise Hamiltonian at any point of $M$.
\end{definition}

Sufficient conditions ensuring that  a canonical action is tubewise Hamiltonian have been given in~\cite{symplectic slice, dual pairs}. For
example, here are two useful results.

\begin{proposition}
\label{iff for strong tubular hamiltonian}
Let $(M,\omega)$ be a symplectic manifold and let $G$ a Lie group with Lie algebra
$\mathfrak{g}$ acting properly and canonically on $M$. For $m\in M$ let
$Y_r:=G \times_{G_{m} } \left(\mathfrak{m}^\ast_r\times V _r\right)$ be the slice
model around the orbit $G\cdot m $ introduced in Proposition~{\rm \ref{design of
symplectic tube at a point}}. If the $G$-equivariant $\mathfrak{g}^\ast$-valued one form
$\gamma
\in \Omega^1(G; \mathfrak{g}^\ast)$ defined by
\begin{equation}
\label{condition for tube wise}
\langle \gamma(g)\left(T _eL _g (\eta) \right), \xi\rangle
:=-\omega(m)\left(\left({\rm Ad } _{g ^{-1}}\xi \right)_M(m),
\eta_M(m)\right) 
\end{equation}
for any $g \in G$ and $\xi, \eta \in \mathfrak{g}$ is exact, then the $G$-action
on $Y _r$ given by $g\cdot [h, \eta, v]:=[gh, \eta,v]$, for any $g \in G $ and any $[h,
\eta, v] \in Y _r$, has an associated standard  momentum map and thus
the $G$-action on $(M, \omega)$ is tubewise Hamiltonian at $m$.
\end{proposition}

\begin{corollary}
\label{conditions that produce strong tubular hamiltonians}
Let $(M,\omega)$ be a symplectic manifold and let $G$ be a 
Lie group with Lie algebra $\mathfrak{g}$ acting properly and
canonically on $M$. If either
\begin{description}
\item[(i)] $H ^1(G)=0 $, or
\item[(ii)] the orbit $G\cdot m$ is isotropic
\end{description}
then the $G$-action on $(M , \omega) $ is tubewise Hamiltonian at $m$.
\end{corollary}

The following result is the analog of Proposition~\ref{local expression j u} in the
non-Abelian setup.

\begin{proposition}
\label{local expression j u non abelian}
Let $(M,\omega)$ be a connected paracompact symplectic manifold and let $G$ be a
compact connected Lie group acting canonically on $M$ in a tubewise Hamiltonian fashion with
closed Hamiltonian holonomy
${\mathcal H}$. Let $\mathbf{K}:M \rightarrow \mathfrak{g}^\ast/ {\mathcal H}$ be a cylinder
valued momentum map for this action and $m  \in M $ arbitrary such that $\mathbf{K}(m)=[
\mu] \in  \mathfrak{g}^\ast/ {\mathcal H} $.  Then, there exists an open neighborhood $U$ of 
$m$ in $M$ and  open neighborhoods $W $ and $V$ of $\mu\in \mathfrak{g}^\ast$ and  $[\mu]\in
\mathfrak{g}^\ast / {\mathcal H}$, respectively,
such that  $\mathbf{K} (U) \subset  V $, $\left.\pi _C\right|_{W}:W
\rightarrow V $ is a diffeomorphism, and 
\begin{equation}
\label{local expression momentum non-abelian}
\left.\pi _C\right|_{W}^{-1}\circ  \mathbf{K}|_{U}= \mathbf{J} _U+c,
\end{equation}
where $c \in \mathfrak{g}^\ast$ is a constant and  $\mathbf{J}_U:U \rightarrow  \mathfrak{g}^\ast$ is a map that in symplectic slice coordinates around the point $m$ has the
expression
\begin{equation}
\label{local expression J non-abelian}
\mathbf{J}_U([g, \rho, v])=\mbox{\rm Ad}^\ast _{g ^{-1}}(\nu+\rho+ \mathbf{J}_V (v)),
\end{equation}
with $\nu \in \mathfrak{g}^\ast$ a constant.
\end{proposition}  

\noindent\textbf{Proof.\ \ } Due to the closedness hypothesis on the Hamiltonian holonomy
${\mathcal H} $, the projection
$\pi_C$ is a local diffeomorphism (see Proposition~\ref{properties of pic closed h}) and
hence there exists an open neighborhood $V$ of  $\mathbf{K} (m)$ in $\mathfrak{g}^\ast/
{\mathcal H} $ and a neighborhood $W$ of some element in the fiber $\pi _C
^{-1}(\mathbf{K} (m) ) \subset \mathfrak{g}^\ast$ such that
$\left.\pi _C\right|_{W}:W
\rightarrow V $ is a diffeomorphism. Let $U$ be the connected component containing $m$ of
the intersection of $\mathbf{K}^{-1}(V) $ with the domain of a symplectic slice chart around
$m$. Given that the $G$-action is by hypothesis tubewise Hamiltonian, the symplectic slice
chart can be chosen so that the restriction of the $G$-action to that chart has a standard 
associated momentum map $ \mathbf{J}_{Y _r}$ that in slice coordinates has, by the
Marle-Guillemin-Sternberg normal form~\cite{marle, GuSt1984},  the expression 
\begin{equation*}
\mathbf{J}_{Y _r}([g, \rho, v])=\mbox{\rm Ad}^\ast _{g ^{-1}}(\nu+\rho+ \mathbf{J}_V (v))+
\sigma(g),
\end{equation*}
with $\nu\in \mathfrak{g}^\ast$ a constant and $\sigma : G \rightarrow \mathfrak{g}^\ast$
the non-equivariance one-cocycle of $\mathbf{J}_{Y _r}$. Since the group $G$ is compact, $\mathbf{J}_{Y _r}$ can be chosen  equivariant and
hence with trivial non-equivariance cocycle $\sigma$ (see~\cite{mre} for the original source of this result, or \cite{hsr}, Proposition 4.5.19).  Let $\mathbf{J}_U $ be the restriction
of that equivariant  momentum map to $U$.
By the definition of the standard momentum map we have that  for any $z \in U $ and any
$\xi\in
\mathfrak{g} $, the map
$\mathbf{J}_U ^\xi:=\langle \mathbf{J} _U, \xi\rangle $ satisfies 
$X_{\mathbf{J}_U ^\xi}(z)= \xi_M (z)$,
with $X_{\mathbf{J}_U ^\xi} $ the Hamiltonian vector field associated to the function
$\mathbf{J}_U ^\xi \in  C^\infty(M) $. With this in mind, it suffices to mimic the proof of
Proposition~\ref{local expression j u} starting from expression~(\ref{local expression
momentum infinitesimal}) to establish the statement of the proposition.
\quad
$\blacksquare$

\begin{theorem}
\label{non abelian convexity}
Let $(M,\omega)$ be a connected paracompact symplectic manifold and let $G$ a
compact connected Lie group acting canonically on $M$ in a tubewise Hamiltonian fashion with
closed Hamiltonian holonomy
${\mathcal H}$. Let $\mathbf{K}:M \rightarrow \mathfrak{g}^\ast/ {\mathcal H}$ be a cylinder
valued momentum map for this action. If $\mathbf{K}$ is a closed map then the intersection
of the image
$\mathbf{K}(M) \subset \mathfrak{g}^\ast/{\mathcal H} $ with $\mathfrak{t}^\ast _+/
{\mathcal H} $  is a weakly convex subset of
$\mathfrak{g}^\ast/ {\mathcal H} $. We think of $\mathfrak{g}^\ast/{\mathcal H}$ and
$\mathfrak{t}^\ast/{\mathcal H}$ as length metric spaces with the length metric naturally
inherited from $\mathfrak{g}^\ast$ {\rm (}see Proposition {\rm \ref{properties of pic closed h})}. If, in
addition, $\mathfrak{t}^\ast_+/ {\mathcal H} $ is uniquely geodesic then $\mathbf{K} (M)\cap
(\mathfrak{t}^\ast _+/ {\mathcal H}) $ is convex, $\mathbf{K} $ has connected fibers, and it
is open onto its image. 
\end{theorem}

\noindent\textbf{Proof.\ \ } First of all,
notice that the closedness of the Hamiltonian holonomy implies, by
Proposition~\ref{properties of pic closed h},  that $\mathfrak{g}^\ast/{\mathcal H}$ and
$\mathfrak{t}^\ast_+/{\mathcal H}$  are complete and locally compact length spaces and that
the projections $\pi_C$ and $\pi_C^+ $ are  local isometries. Moreover, the
identification $\mathfrak{t}^\ast
_+/ {\mathcal H}\simeq (\mathfrak{g}^\ast/  {\mathcal H})/G $, introduced in part \textbf{(iii)}  of Proposition~\ref{properties of pic closed h}  and the diagram~(\ref{commuting
pi}) allow us to think of
$\pi _C^+ $ as the restriction of $\pi _C $ to $\mathfrak{t}^\ast _+ $. Consequently, if 
$V \subset \mathfrak{g}^\ast/ {\mathcal H}$ is an open set such that 
$\left.\pi _C\right|_{W}:W
\rightarrow V $ is an isometric diffeomorphism then  $\left.\pi
_C^+\right|_{\pi _G(W)}:\pi_G(W)\simeq W\cap
\mathfrak{t}^\ast _+
\rightarrow \pi _G^+(V) \simeq (W\cap
\mathfrak{t}^\ast _+)/ {\mathcal H}$ is an isometry. Notice that $\pi_G(W)\simeq
W\cap
\mathfrak{t}^\ast _+ \subset
\mathfrak{t}^\ast _+
$ is an open subset of $\mathfrak{t}^\ast _+ $  since $\pi_G  $ is an open map.

Using the identification $\mathfrak{t}^\ast
_+/ {\mathcal H}\simeq (\mathfrak{g}^\ast/  {\mathcal H})/G $ we can
study the convexity properties of the intersection $\mathbf{K} (M)\cap
(\mathfrak{t}^\ast _+/ {\mathcal H}) $ by looking at the convexity properties of the image of
the map $k:= \pi _G^+ \circ \mathbf{K} :M \rightarrow  \mathfrak{t}^\ast  _+/ {\mathcal H}$.
We will do so by  applying the
Local-to-Global
Principle for length spaces (Theorem~\ref{lokal to global metric}) to $k$, that
is, by showing that 
$k$ is locally open onto its image, locally fiber connected, and has local
convexity data. 
By Proposition~\ref{local expression j u non abelian}
there exists an open neighborhood $U$ of 
$m$ in $M$ and an open neighborhood $V$ of  $[\mu]$ in $\mathfrak{g}^\ast / {\mathcal H}$
such that  $\mathbf{K} (U) \subset  V $, $\left.\pi _C\right|_{W}:W
\rightarrow V $ is a diffeomorphism, and 
$
\left.\pi _C\right|_{W}^{-1}\circ  \mathbf{K}|_{U}= \mathbf{J} _U+c,
$
with $c \in \mathfrak{g}^\ast$ a constant and  $\mathbf{J}_U:U \rightarrow 
\mathfrak{g}^\ast$ a map that  has the
expression~(\ref{local expression J non-abelian}). If we apply $\pi_G $ to both sides of
this equality and we use the commutativity of diagram~(\ref{commuting pi}) and the remarks
above we obtain that
\[
\left(\left.\pi
_C^+\right|_{\pi _G(W)}\right)^{-1}\circ  k|_{U}= j _U+ \pi_G(c),
\]
with $j _U:= \pi _G \circ \mathbf{J}_U$. Two results due to Sjamaar~\cite[Theorem
6.5]{sjamaar} and Knop~\cite[Theorem 5.1]{knop} show that $j _U $ is locally open onto its
image, locally fiber connected, and has local convexity data. Consequently, since $\left.\pi
_C^+\right|_{\pi _G(W)} $ is an isometry, Lemma~\ref{convex image} guarantees
that $k| _U $ also shares those three local properties. The statement of the theorem follows
then as a consequence of Theorem~\ref{lokal to global metric}. \quad $\blacksquare$

\begin{remark}
\normalfont
The classical convexity theorem of Kirwan~\cite{kirwan convexity}  states that if $G$ is a
compact connected Lie group acting canonically on the compact connected symplectic
manifold $(M, \omega )$ and this action has an associated standard  coadjoint equivariant momentum map
$\mathbf{J}:M \rightarrow  \mathfrak{g}^\ast$, then $\mathbf{J}(M)\cap 
\mathfrak{t}^\ast _+ $  is a compact convex polytope. The convexity part of this theorem 
can be obtained from Theorem~\ref{non abelian convexity} since the global existence of a
standard momentum map implies that the action is, in particular, tubewise Hamiltonian and that
a cylinder valued momentum map for it is the momentum map  $\mathbf{J}:M \rightarrow 
\mathfrak{g}^\ast$ (${\mathcal H} =\{0\} $ in this case). Since the manifold $M$ is by
hypothesis compact then
$\mathbf{J} $ is necessarily a closed map and, moreover, in this case $\mathfrak{t}^\ast _+/
{\mathcal H}= \mathfrak{t}^\ast _+ $ is uniquely geodesic. Consequently, by Theorem~\ref{non abelian convexity}, $\mathbf{J}(M)\cap 
\mathfrak{t}^\ast _+ $ is convex, $\mathbf{J} $  has connected fibers, and it is a $G$-open
map onto its image. See also~\cite{polytope paper} where the same result was obtained in a different manner. 
\end{remark}

\section{Appendix: metric and length spaces}
\label{Appendix: metric and length spaces}

In this appendix we collect the standard results and we fix the notations that we use when dealing with metric and length spaces.  Most of the quoted statements below can be
found in Bridson and Haefliger \cite{bridson} and Burago et al. 
\cite{burago}.

\begin{definition}\label{metric space}
Let $X$ be an arbitrary set. A function $d: X \times X \to \mathbb
R \cup \{\infty\}$ is called a \textbf{metric} on $X$ if the following
conditions are satisfied for all $x, y, z \in X$:
\begin{description}
  \item[(i)]  Positiveness: $d(x, y) >0$ if $x \not = y$ and $d(x,x)=0$;
  \item[(ii)]  Symmetry: $d(x, y) = d(y,x)$;
  \item[(iii)]  Triangle inequality: $d(x, z) \leq d(x, y)+d(y,z)$.
\end{description}
A \textbf{metric space} is a pair $(X, d)$ with $d$ a metric on the
set
$X$.
\end{definition}

Let $(X, d)$ and $(X^\prime, d^\prime)$ be two metric spaces. A map
$f : X \to X^\prime$ is  called \textbf{distance-preserving} if
$d^\prime (f(x), f(y))=d(x,y)$ for all $x, y \in X$. A bijective
distance-preserving map is called an \textbf{isometry}. Note that a
distance-preserving map is always injective.

Given $x \in X$ and $r>0$ the \textbf{open ball} of radius
$r$ and center $x$ is defined to be the set $B(x, r): = \{y \in X  \ | \ d(x, y)
<r\}$. Similarly, the \textbf{closed ball} of radius $r$ and center $x$ is the
set $\overline B (x, r) : = \{y \in X \ | \ d(x,y) \leq r\}$. On
$(X,d)$ the topology given by the metric $d$ is, by definition, the topology
whose basis of neighborhoods at every point $x \in X $ is the
collection of all open balls $B(x, r)$ for $r >0$. Thus, a set
$U$ in the metric space $(X, d)$ is open if and only if for every
point $x \in U$ there exists $r>0$ such that
$B(x, r) \subset U$. It is easy to see that $\overline B(x, r)$ are closed
sets in the metric topology but, in general, they are strictly larger
that the closure $\overline{B(x,r)}$ of the open balls in the same
topology.

\begin{definition}\label{complete metric space}
Let $(X, d)$ be a metric space. A sequence $\{x_n\}$ is called
\textbf{Cauchy sequence} if $d(x_n, x_m) \to 0$ as $n, m \to \infty$. A
metric space is called \textbf{complete} if  every Cauchy sequence has a
convergent subsequence.
\end{definition}

Let $(X,d)$ be a metric space and $Y \subset X$ a subset. Then the
restriction $d_Y$ of the metric $d$ to $Y \times Y$ defines a metric on $Y$. Thus, $(Y, d_Y)$ is a metric space whose metric topology coincides with the relative topology induced from the metric topology of $X$. Also, if $Y$ is a complete metric space relative to the induced
metric $d_Y$, then $Y$ is closed in $X$ and if $(X,d)$ is a complete metric space and $Y$ is closed in $X$, then $Y$ is complete.

\medskip

Let $(X,d)$ be a metric space. A \textbf{curve} or a \textbf{path} in
$X$ is a continuous map $c: I\to X$ with $I$ a connected interval
of $\mathbb R$. If $c_1:[a_1, b_1]\to X$ and $c_2:[a_2,b_2]\to X$
are two paths such that $c_1(b_1)=c_2(a_2)$, their \textbf{concatenation} is
the path $c:[a_1, b_1+b_2-a_2]\to X$ defined by $c(t)=c_1(t)$ if
$t \in [a_1, b_1]$ and $c(t) = c_2 (t+a_2-b_1)$ if $t \in [b_1,
b_1+b_2-a_2]$.

\begin{definition}\label{d5}
The \textbf{length} $l_d (c)$ of a curve $c:[a, b]\to X$ induced by the
metric $d$ is
\begin{equation*}
    l_d(c): = \sup_{\Delta_n}
    \sum_{i=0}^{n-1} d(c(t_i), c(t_{i+1})),
\end{equation*} where the supremum is taken over all possible
partitions $\Delta_n:a=t_0\leq t_1 \leq \dots \leq t_n=b$ of the
interval $[a,b] \subset \mathbb R$.
\end{definition}

Consequently, the length of a curve is a non-negative number or it is
infinite. The curve $c$ is said to be \textbf{rectifiable} if its
length is finite. 

Next, we recall several properties of the
length of a curve in a metric space.

\smallskip

\begin{description}
  \item[(i)]  $l_d (c) \geq d(c(a), d(c(b)))$, for any path $c:[a,b]\to
  X$.
  \item[(ii)]  If $\phi:[a^\prime, b^\prime]\to [a,b]$ is an onto
monotone
  map, then $l_d(c)=l_d(c\circ \phi)$.
  \item[(iii)]  Additivity: if $c$ is the concatenation of two paths $c_1$
  and $c_2$ then $l_d (c) = l_d (c_1) +l_d (c_2)$.
  \item [(iv)] If $c$ is rectifiable of length $l$, then the function
  $\lambda:[ a,b]\to [0,l]$ defined by $\lambda(t) = l_d (c_{[a,t]})$
is a continuous weakly monotone function.
  \item[(v)]  Reparametrization by arc length: if $c$ and $\lambda$ are as in
  as in the previous point, then there is a unique path
$\widetilde{c}:[0,l]\to X$ such that
  \begin{equation*}
    \widetilde{c} \circ \lambda = c \quad \mbox{and} \quad
l_d(\widetilde{c}_{[0,t]})=t.
  \end{equation*}
	\item [(vi)] Lower semicontinuity: let $(c_n)$ be a sequence of paths
  $[a,b] \to X$ converging uniformly to a path $c$. If $c$
  rectifiable, then for every $\varepsilon>0$, there exists an integer $N_\varepsilon$
  such that
  \begin{equation*}
    l_d(c) \leq l_d(c_n) +\varepsilon
  \end{equation*} whenever $n>N_\varepsilon$.
\end{description}

\smallskip

Next we will introduce the notion of length metric or inner
metric. It is well known that every Riemannian metric on a manifold
induces a length. Unfortunately, Riemannian metrics can be defined only in the
differentiable setting. As we will see,  length metrics  share many properties with
Riemannian  metrics but they can be  defined in more general settings.

\begin{definition}\label{length space}
Let $(X,d)$ be a metric space. The distance $d$ is said to be a
\textbf{length metric} or an \textbf{inner metric} if the distance
between every pair of points $x, y \in X$ is equal to the infimum
of the length of rectifiable curves joining them. If there are no
such curves then, by definition, $d(x,y) =\infty$. If $d$ is a length
metric then
$(X,d)$ is called a \textbf{length space} or an \textbf{inner space}.
\end{definition}

Other authors refer to length spaces as \textbf{path metric spaces} (see, e.g., Gromov \cite{gromov}). For various properties and characterizations
of length metrics see Bridson and Haefliger
\cite{bridson}, Burago et al. \cite{burago}, and Gromov
\cite{gromov}.

Having a metric space $(X,d)$ we can always construct a length
metric $\overline d$ induced by the initial metric $d$ in the
following way,
\begin{equation*}
    \overline d (x, y) : = \inf_{\emph{R}_{x,y}}l_d(\gamma),
\end{equation*} 
where $R_{x,y} := \{\text{all rectifiable curves
connecting } x \text{ and } y\}$. If there are no such curves then
we set $\overline d(x,y)=\infty$.

\begin{proposition}[Bridson and Haefliger \cite{bridson}]
The induced length metric has the following properties:
\begin{description}
  \item [(i)] $\overline d$ is a metric.
  \item [(ii)] $\overline d(x,y) \geq d(x,y)$ for all $x,y \in X$.
  \item [(iii)] If $c:[a,b]\to X$ is continuous with respect to the topology
  induced by $\overline{d}$, then it is continuous with respect to the
  topology induced by $d$. (The converse is false, in general).
  \item [(iv)] If a map $c:[a,b]\to X$ is a rectifiable curve in $(X,d)$,
  then it is a continuous and rectifiable curve in $(X, \overline d)$.
  \item [(v)] The length of a curve $c:[a,b]\to X$ in $(X, \overline d)$ is the
  same as its length in $(X,d)$.
  \item [(vi)] $\overline{\overline d}=\overline{d}$.
\end{description}
\end{proposition}

The assertion in point \textbf{(iii)} of the above proposition is a consequence of the fact that the topology induced by the metric $d$ is coarser
than the topology induced by the metric $\overline d$. Note that
$(X,d)$ \textit{is a length space if and only if\/} $\overline d =d$.
\medskip

Classical examples of a length spaces are Riemannian manifolds.
Let $(X,g)$ be a Riemannian manifold and $c:[a,b] \to X$ a piecewise
differentiable path. The Riemannian length $l_g(c)$ is defined as
\begin{equation*}
    l_g(c):= \int_a^b\sqrt{g_{ij}(t)\dot c^i(t)\dot c^j(t)}.
\end{equation*}

\begin{proposition}\label{Riemannian geometry}
Let $X$ be a connected Riemannian manifold. Given $x,y \in X$, let
$d(x,y)$ be the infimum of the Riemannian length of piecewise
continuously differentiable paths $c:[0,1]\to X$ such  that
$c(0)=x$ and $c(1)=y$. Then
\begin{description}
  \item [(i)] $d$ is a metric on $X$.
  \item [(ii)] The topology on $X$ defined by this distance is the same as
  the given manifold topology on $X$.
  \item [(iii)] $(X,d)$ is a length space.
\end{description}
\end{proposition}

\begin{definition}\label{shortest path}
A curve $c:[a,b]\to (X,d)$ is called a \textbf{shortest path} if its
length is minimal among all the curves with the same endpoints. 
Shortest paths in length spaces are also called
\textbf{distance minimizers}.
\end{definition}

If $(X,d)$ is a length space then a curve $c:[a,b]\to X$ is a
shortest path if and only if its length is equal with the distance
between endpoints, that is,
$l_d(c)=d(c(a), c(b))$.
\medskip 

Next we will introduce the notion of geodesic in length spaces
that generalizes the one in Riemannian geometry.

\begin{definition}\label{geodesic}
Let $(X,d)$ be a length space. A curve $c:I \subset \mathbb R \to
X$ is called \textbf{geodesic} if for every $t \in I$ there exist a
subinterval $J$ containing a neighborhood of $t$ in  $I$ such that
$c|_{J}$ is a shortest path. In other words, a geodesic is a curve
which is locally a distance minimizer. A length space  $(X,d)$ is
called a \textbf{geodesic metric space} if for any two points $x,y
\in X$ there exists a shortest path between $x$ and $y$.
\end{definition}

Clearly in a length space a shortest path is a geodesic. The
extension of the Hopf-Rinow theorem from Riemannian geometry to the
case of length metric spaces is due Cohn-Vossen and is a
key result in this paper. Its proof can be found in
Bridson and Haefliger \cite{bridson} or Burago et al. \cite{burago}.

\begin{theorem}\label{Hopf-Rinow}
(Hopf-Rinow-Cohn-Vossen) For a locally compact length space
$(X,d)$, the following assertions are equivalent:
\begin{description}
  \item [(i)] $X$ is complete,
  \item [(ii)] every closed metric ball in $X$ is compact.
\end{description}
If one of the above assertions holds, then for any two points
$x,y \in X$ there exists a shortest path connecting them. In other
words, $(X, d)$ is a geodesic metric space.
\end{theorem}

\begin{corollary}\label{c1}
Every complete, connected, Riemannian manifold is a geodesic
metric space.
\end{corollary}

Having introduced the notion of shortest path one can  define the
key concept of metric convexity.

\begin{definition}\label{metric convex}
A subset $C$ in a metric space $(X, d)$ is said to be \textbf{convex}
if the restriction of $d$ to $C$ is a finite length metric.
\end{definition}

For the case of geodesic metric spaces we have the following
characterization of convex sets.

\begin{lemma}\label{caracterizare}
Let $(X,d)$ be a geodesic metric space. Then a subset $C \subset X$ is
convex if and only if for any two points $x,y\in C$ there exists a
rectifiable shortest path $\gamma$ connecting $x$ and $y$ which is
entirely contained in $C$.
\end{lemma}

\noindent\textbf{Proof.\ \ } Let $d_C$ be the restriction of $d$ to the subset $C$
and assume $C$ is a convex subset of $(X, d)$. Since $d_C = \overline
{d_C}$ and $d_C<\infty$ we have the following equality

\begin{equation*}
\inf_{\emph{R}_{x,y}} l_d (\gamma)=d_C(x,y)=\overline{d_C}(x,y)=
\inf_{\emph{R}^C_{x,y}}l_d (\gamma),
\end{equation*}
where $\emph{R}_{x,y}^C:=\{\gamma :[a,b]\to C \ | \ \mbox{$\gamma$
rectifiable curve entirely contained in}\; C \; \mbox{with}\;
\gamma (a) =x \; \mbox{and} \; \gamma (b) =y\}.$ Observe that
$\emph{R}_{x,y}^C\neq \varnothing$ because $C$ is a length space
whose metric $d_C$ is finite. This
shows that taking the infimum over $\emph{R}^C_{x,y}$ completely
determines $d_C(x,y)$. Because $(X,d)$ is a geodesic metric space
this infimum is attained for a curve entirely contained in $C$.

Conversely, suppose that for any two points $x,y\in C$ there
exists a rectifiable shortest path between $x$ and $y$ which
is entirely contained in $C$. We have to prove that $d_C<\infty$ and that 
$d_C=\overline{d_C}$. To prove the first claim note that, by the working hypothesis, there exists a rectifiable path $\gamma_{0}$
entirely belonging to $C$ connecting $x$ and $y$ such that
${d_C}(x,y)=l_{d}(\gamma_{0})<\infty$. Since
$\inf_{\emph{R}_{x,y}} l_d (\gamma)$ is attained for paths entirely
belonging to $C$ we obtain the equality $d_C=\overline{d_C}$, which proves the second claim.
\quad $\blacksquare$

\begin{lemma}\label{convex image}
Let $(X,d_X)$ and $(Y, d_Y)$ be two geodesic metric spaces  and
$f:X\to Y$ a distance-preserving map. If $d_X$ is finite then
the image $\operatorname{Im} (f):=f(X)$ is a convex subset of $Y$.
\end{lemma}

\noindent\textbf{Proof.\ \ } Since $f$ is a distance preserving map it is injective. Let $y_1, y_2 \in f(X)$ and $x_1, x_2 \in X$ be the corresponding
preimages. Then  there exists a shortest
geodesic between $x_1$ and $x_2$, namely a curve $c:[a,b]\to X$ satisfying $c(a) = x_1$,
$c(b)=x_2$ and
$l_{d_X}(c)=d_X(x_1,x_2)<\infty$. Since $f$ is distance preserving
we have that $l_{d_X}(c) = l_{d_Y} (f\circ c)$. Consequently
$l_{d_Y} (f\circ c)=d_{Y}(y_1, y_2)<\infty$ which proves that
$f\circ c$ is a shortest geodesic connecting $y_1$ and $y_2$
entirely contained in $f(X)$. Consequently, $f(X)$ is a convex
subset of $Y$.
\quad $\blacksquare$

\medskip

\noindent\textbf{Acknowledgments}   
We thank A. Alekseev for patiently answering our questions about Lie group valued momentum maps over the years.
P. B. has been supported by a grant of the R\'egion de Franche-Comt\'e (Convention 051004-02) during his stay at the Universit\'e de Franche-Comt\'e, Besan\c{c}on, which made possible this collaboration. T.S.R. was partially supported by a Swiss National Science Foundation grant.

\noindent {\sc P. Birtea} \\
Departamentul de Matematic\u a, Universitatea de Vest,
RO--1900 Timi\c soara, Romania.\\
Email: {\sf birtea@math.uvt.ro}
\medskip

\noindent {\sc J.-P. Ortega} \\
D\'epartement de Math\'ematiques de Besan\c{c}on,
Universit\'e de Franche-Comt\'e, UFR des Sciences et
Techniques, 16 route de Gray, F--25030 Besan\c{c}on c\'edex,
France.\\
Email: Juan-Pablo.Ortega@math.univ-fcomte.fr
\medskip

\noindent {\sc T.S. Ratiu}\\
Section de Math\'ematiques and Bernoulli Center,
{\'E}cole Polytechnique F{\'e}d{\'e}rale de Lausanne,
CH--1015 Lausanne,
Switzerland.\\
Email: {\sf tudor.ratiu@epfl.ch}

\end{document}